\documentclass[11pt]{article}

\usepackage{amsfonts}

\usepackage{CJK}
\usepackage{amsfonts,amssymb,amsmath,mathrsfs,multirow,booktabs}
\usepackage{xcolor,soul}
\usepackage{color,latexsym,amsfonts}
 \newcommand{\qed}{\hfill\rule{2mm}{3mm}\vspace{4mm}}

 \newtheorem{theorem}{Theorem}[section]
 \newtheorem{lemma}[theorem]{Lemma}
 \newtheorem{corollary}[theorem]{Corollary}
 \newtheorem{proposition}[theorem]{Proposition}
 \newtheorem{example}[theorem]{Example}
 \newtheorem{Definition}[theorem]{Definition}
 \newtheorem{remark}[theorem]{Remark}
 \newtheorem{condition}[theorem]{Condition}
 \newtheorem{conjecture}[theorem]{Conjecture}

 \def\blemma{\begin{lemma}\sl{}\def\elemma{\end{lemma}}}
 \def\btheorem{\begin{theorem}\sl{}\def\etheorem{\end{theorem}}}
 
 \def\bdefinition{\begin{Definition}\sl{}\def\edefinition{\end{Definition}}}
 \def\bproposition{\begin{proposition}\sl{}\def\eproposition{\end{proposition}}}
 \def\bremark{\begin{remark}\sl{}\def\eremark{\end{remark}}}
 \def\bcondition{\begin{condition}\sl{}\def\econdition{\end{condition}}}

 \def\beqlb{\begin{eqnarray}}\def\eeqlb{\end{eqnarray}}
 \def\beqnn{\begin{eqnarray*}}\def\eeqnn{\end{eqnarray*}}

 \def\mbb{\mathbb}\def\mbf{\mathbf}

 \def\<{\langle}\def\>{\rangle}

 \def\ar{&\!\!}

 \def\eqref#1{{\rm(\ref{#1})}}

 \def\proof{\noindent{\it
 Proof.~}}\def\qed{\hfill$\Box$\medskip}

\def\e{{\mbox{\rm e}}}
\def\<{\left<}\def\>{\right>}

 \def\mbb{\mathbb}
  \def\mbf{\mathbf}
\newcommand{\dd}{\mathrm{d}}

\newfam\msbmfam\font\tenmsbm=msbm10\textfont
\msbmfam=\tenmsbm\font\sevenmsbm=msbm7
\scriptfont\msbmfam=\sevenmsbm

\def\<{\left<}\def\>{\right>}
\def\({\left(}\def\){\right)}

\begin{document}

\centerline{\Large\bf  Behaviors near infinity for mutually}
\medskip
\centerline{\Large\bf enhancing continuous-state population dynamics}

\bigskip

\centerline{
Jie Xiong\footnote{Department of Mathematics and SUSTech International center for Mathematics, Southern University of Science \& Technology, Shenzhen, China.
Supported by National Key R\&D Program of China (No.~2022YFA1006102) and
National Natural Science Foundation of China Grant 12471418.  Email: xiongj@sustech.edu.cn},
Xu Yang\footnote{School of Mathematics and Information Science,
North Minzu University, Yinchuan, China. Supported by
NSFC (No.~12471135). Email: xuyang@mail.bnu.edu.cn. Corresponding author.} and Xiaowen Zhou\footnote{Department of Mathematics and Statistics, Concordia
University, Montreal, Canada.
Supported by NSERC (RGPIN-2021-04100).
Email: xiaowen.zhou@concordia.ca.}}

\bigskip\bigskip

{\narrower{\narrower

\noindent{\bf Abstract.}
In this paper we consider a two-dimensional generalized continuous-state branching process with  mutually enhancing two-way interactions
characterized by two stochastic differential equations driven by Brownian motions and spectrally positive $\alpha$-stable random measures.
Such continuous-state population dynamics can also be identified as a stochastic Lotka-Volterra type population model.
Infinite behavior such as
explosion/nonexplosion for this population
and
staying infinite/coming down from infinity are derived under various
conditions on the coefficients involved in the model and
the conclusions are rather sharp.

\medskip

\textit{Mathematics Subject Classifications (2010)}: 60J80; 92D25; 60G57; 60G17.

\medskip

\textit{Key words and phrases}:
continuous-state branching process,
mutually enhancing interaction,
stochastic Lotka-Volterra type population,
staying infinite, coming down from infinity,
explosion, nonexplosion.
\par}\par}

\section{Introduction and main results}\label{intro}

\setcounter{equation}{0}

Continuous-state branching process (CSBP for short),
which arises as time-space scaling limit of
discrete Bienaym\'e-Galton-Watson (BGW for short) processes
describing the evolutionary behavior of the population,
is a nonnegative real-valued Markov process.
 It can also be obtained as time-changed of a L\'evy process with no negative jumps.
Stochastic differential equation (SDE for short) theory plays an important role in the study of CSBPs.
We refer to Kyprianou \cite{Ky} and Li \cite{Li2020} for nice introduction on CSBPs and the associated SDEs.

The classical BGW process needs the species satisfying the additive branching property.
However, in general population systems, due to species competition and limited resources, species no longer meet the additive branching requirements.
Thus generalized CSBPs such as competition and nonlinear branching
have been proposed in recent years.
In particular, \cite{Pardoux16} introduces CSBPs with competition
which is characterized as a SDE called Dawson-Li equation.
The CSBPs with polynomial branching are considered in \cite{LiP19},
and with competition and more general branching
characterized as SDEs driven by Gaussian white noises and compensated measures
are discussed in \cite{LYZh}.

A classical question in the theoretical study of Markov processes is
long-term behaviour such as explosion, extinction and coming down from infinity
and so on. The phenomenon of coming down from infinity has been extensively studied for coalescing Brownian motions; see, e.g., \cite{BarnesMytnikSun}.
For classical CSBP, the explosion/nonexplosion conditions are given by
an integral test on Laplace exponent of the associated branching mechanism
in \cite{Grey74}. The explosion/nonexplosion and staying infinite/coming down from infinity
conditions for the CSBPs with polynomial branching are studied
by Lamperti type transformations in \cite{LiP19},
and for the general CSBPs are given by martingale approaches in \cite{LYZh}.
The speed of coming down from infinity and explosion  for the CSBPs with polynomial branching are studied by \cite{FoucartLiZhou21} and \cite{LiZh21},
respectively.
Recently, the explosion condition for more general time-inhomogeneous SDE with jumps
is discussed in \cite{ChenYangZhou}.

Compared with one-dimensional generalized CSBPs,
it is interesting and challenging to study two-dimensional interacting
generalized CSBPs
characterized by the following SDEs system:
	\begin{equation}\label{1.1a}
		\left\{
		\begin{aligned}
			X_t &=X_0 +\int_0^t\Theta_1(X_s,Y_s)\dd s+\int_0^t\gamma_{10}(X_s)\dd s+ \int_0^t\sqrt{\gamma_{11}(X_s)}\dd B_1(s) \\
			&~~~~
			+\int_0^t\int_0^\infty \int_0^{\gamma_{12}(X_{s-})}z\tilde{N}_1(\dd s,\dd z,\dd u),\\
			Y_t &=Y_0+\int_0^t \Theta_2(Y_s,X_s)\dd s
			+\int_0^t\gamma_{20}(Y_s)\dd s +\int_0^t\sqrt{\gamma_{21}(Y_s)}\dd B_2(s) \\
			&~~~~
			+\int_0^t\int_0^\infty \int_0^{\gamma_{22}(Y_{s-})}z\tilde{N}_2(\dd s,\dd z,\dd u),
		\end{aligned}
		\right.
	\end{equation}
where $(B_1(t))_{t\ge0}$ and $(B_2(t))_{t\ge0}$ are Brownian motions,
and $\tilde{N}_1(\dd s,\dd z,\dd u)$ and $\tilde{N}_2(\dd s,\dd z,\dd u)$
are compensated Poisson random measures.
If the coefficients $\Theta_i,\gamma_{ij}$ are linear functions,
the solution to SDEs system \eqref{1.1a} is
a two-dimensional CSBP
and studied in \cite{Watanabe,Cattiaux,Ma2013}.
SDEs system \eqref{1.1a} is also
a stochastic Lotka-Volterra type population dynamical system.
If $\gamma_{i1}(x)=b_{i1}x^2$ and $\gamma_{i2}=0$ for $x>0$ and $i=1,2$
in \eqref{1.1a},
SDEs system \eqref{1.1a} is called the
Lotka-Volterra model in random environments and studied in
 \cite{Bao,Evans,BaoShao,DuDangYin} and the references.
If $\gamma_{10}(x),\gamma_{20}(y)\le0$,
then rather sharp conditions on extinction-extinguishing dichotomy is
given  in \cite{RXYZ19}
for the functions $\Theta_1=0$ and $\Theta_2<0$,
and in \cite{XYZh24} for general mutual competed case $\Theta_1,\Theta_2<0$.
Recently, \cite{XYZh25} discusses the extinction/non-extinction conditions
for mutual enhancing case $\Theta_1,\Theta_2>0$.
While the research on infinite behavior such as explosion/nonexplosion and
staying infinite/coming down from infinity for SDEs system \eqref{1.1a}
with mutual enhancing is more challenging and the conclusions are fewer
as far as we know.
In view of this, {\it our purpose} for this paper is
to establish the conditions on
explosion/nonexplosion and
staying infinite/coming down from infinity.

For simplicity and readability in this paper we only study the special form of \eqref{1.1a} with power coefficients and special densities for $\tilde{N}_i(\dd s,\dd z,\dd u)$, that is the following SDE system:
\begin{equation}\label{1.1}
  \left\{
   \begin{aligned}
   X_t &=X_0 +a_1\int_0^tX_s^{\theta_1}Y_s^{\kappa_1} \dd s-b_{10}\int_0^tX_s^{r_{10}}\dd s+\int_0^t\sqrt{2 b_{11}X_s^{r_{11}}}\dd B_1(s) \\
   &~~~~
+\int_0^t\int_0^\infty \int_0^{b_{12}X_{s-}^{r_{12}}}z\tilde{N}_1(\dd s,\dd z,\dd u),\\
   Y_t &=Y_0+a_2\int_0^t Y_s^{\theta_2}X_s^{\kappa_2}\dd s
   -b_{20}\int_0^tY_s^{r_{20}}\dd s +\int_0^t\sqrt{2 b_{21}Y_s^{r_{21}}}\dd B_2(s) \\
   &~~~~
+\int_0^t\int_0^\infty \int_0^{b_{22}Y_{s-}^{r_{22}}}z\tilde{N}_2(\dd s,\dd z,\dd u),
   \end{aligned}
   \right.
  \end{equation}
where the constants $a_i,\kappa_i>0$ and $\theta_i,r_{ij},b_{ij}\ge0$
satisfying $b_{i1}+b_{i2}>0$ for  $i=1,2$ and $j=0,1,2$.
For $i=1,2$, $(B_i(t))_{t\ge0}$ is a Brownian motion and
$\tilde{N}_i(\dd s,\dd z,\dd u)$ is a compensated Poisson
random measure with intensity $\dd s\mu_i(\dd z)\dd u$.
Here
$\mu_i(\dd z)=\frac{\alpha_i(\alpha_i-1)}{\Gamma(\alpha_i)\Gamma(2-\alpha_i)}z^{-1-\alpha_i}
1_{\{z>0\}}\dd z$
for $\alpha_i\in(1,2)$, $i=1,2$, and $\Gamma$ denotes the Gamma function.
In this paper we assume that
$(B_1(t))_{t\ge0}$, $(B_2(t))_{t\ge0}$, $\{\tilde{N}_1(\dd s,\dd z,\dd u)\}$ and $\{\tilde{N}_2(\dd s,\dd z,\dd u)\}$ are independent of each other.

The main approaches of this paper are to develop the criteria of
explosion/nonexplosion
and staying infinite/coming down from infinity
for two-dimensional process in Section 2 which is an adaption of the method for Chen's criteria
on the uniqueness problem of Markov jump processes.
These criteria are first established in \cite{Chen1986a,Chen1986b} and can also be found in \cite[Theorems 2.25 and 2.27]{Chen04}.
The key in applying these criteria is to find appropriate test functions
for which are the form of two-dimensional power type functions
and often without an obvious intuition.

The main results on the behaviors near infinity for system \eqref{1.1} are given in Theorems \ref{tt1.3a}--\ref{tt1.3}, \ref{t1.4} and \ref{t1.8}--\ref{t1.111}.
All possible outcomes have been covered in special case (see Remark \ref{r1.2}).
The novelty of this paper lies in giving sufficient conditions for both
$X$ and $Y$ to either come down from infinity or stay infinite, which differs from the findings in \cite{XYZh24} and \cite{XYZh25}.

We introduce the following notations.
For any generic stochastic process $V:=(V(t))_{t\ge0}$ and constant $w>0$, let
 \beqnn
\tau_0(V):=\inf\{t\ge0:V(t)=0\},\quad
\tau_w^-(V):=\inf\{t\ge0:V(t)<w\}
 \eeqnn
and
 \beqnn
\tau_w^+(V):=\inf\{t\ge0:V(t)>w\}
 \eeqnn
with the convention $\inf\emptyset=\infty$.
Let $\tau_0:=\tau_0(X)\wedge\tau_0(Y)$,
$\tau_w^-:=\tau_w^-(X)\wedge\tau_w^-(Y)$,
$\tau_w^+:=\tau_w^+(X)\wedge\tau_w^+(Y)$
and $\tau_\infty:=\lim_{n\to\infty}\tau_n^{+}$.
In the following we state the definition of solution
to SDEs system \eqref{1.1}, which is defined before the minimum of
the processes first time of either hitting zero or explosion.

\bdefinition\label{def}
By a solution to SDE \eqref{1.1} we mean that a two-dimensional c\'adl\'ag
process $(X,Y):=((X_t,Y_t))_{t\ge0}$ satisfies
SDE \eqref{1.1} up to $\gamma_n:=\tau_{1/n}^-
\wedge\tau_n^+$ for each $n\ge1$ and
$X_t=\limsup_{n\to\infty}X_{\gamma_n-}$
and $Y_t=\limsup_{n\to\infty}Y_{\gamma_n-}$
for $t\ge\lim_{n\to\infty}\gamma_n$.
\edefinition

The above definition means that the solution is nonnegative
and $0$ and $\infty$ are absorbing states of the system.
The existence and pathwise uniqueness to SDEs system \eqref{1.1} can be gotten by the same arguments as in \cite[Lemma A.1]{RXYZ19}.
In this paper we always assume that the $(X,Y)$ is the unique solution to SDEs system \eqref{1.1}, and consequently, $(X,Y)$ has the strong Markov property.
We also assume that $X_0,Y_0>0$ are deterministic and all the stochastic processes are defined on the same filtered probability space $(\Omega,\mathscr{F},\mathscr{F}_t,\mathbf{P})$.
Let $\mathbf{E}$ denote the corresponding expectation,
$\mathbf{P}_{(x,y)}$ be the law of a
process started at $(x,y)$.
The system $(X,Y)$ {\it stays infinite} if
$\lim_{X_0,Y_0\to\infty} \mbf{P}_{(X_0,Y_0)}\{\tau_a^-(X)\wedge\tau_a^-(Y) <t\}=0$ for all $a,t>0$.
The system $(X,Y)$ {\it comes down from infinity} if
$\lim_{a\to\infty}\lim_{X_0,Y_0\to\infty} \mbf{P}_{(X_0,Y_0)}\{\tau_a^-(X)\vee \tau_a^-(Y)<t\}=1$
and all $t>0$.
For $i=1,2$ set
 \beqnn
 r_i
:=
\max_{j=0,1,2}\big\{r_{ij}-\varrho_{ij}:b_{ij}\neq0\big\},\qquad
 b_i
:=
\sum_{j=0,1,2}b_{ij}1_{\{ r_i=r_{ij}-\varrho_{ij}\}}.
 \eeqnn
where $\varrho_{ij}:=j+1$ for $j=0,1$
and $\varrho_{i2}:=\alpha_i$ for all $i=1,2$.

By \cite[pp.2536--2537]{LYZh}, the one-dimensional nonlinear CSBP
(with $a_1X_s^{\theta_1}Y_s^{\kappa_1}$ replaced by $a_1X_s^{\theta_1}$ in \eqref{1.1}) is explosion if $r_1\vee0<\theta_1-1$
and nonexplosion if $r_1\vee0>\theta_1-1$,
and stays infinite if $r_1<(\theta_1-1)\vee0$
and comes down from infinity if $r_1>(\theta_1-1)\vee0$.
The following results present that the explosion
and staying infinite are consistent with
the corresponding one-dimensional CSBP excluding the critical case, which is obvious intuitively
and means that interaction terms
$a_1X_s^{\theta_1}Y_s^{\kappa_1}$ and $a_2Y_s^{\theta_2}X_s^{\kappa_2}$
are useless.
\btheorem\label{tt1.3a}
\begin{itemize}
\item[{\normalfont(i)}]
$\lim_{X_0,Y_0\to\infty}\mbf{P}_{(X_0,Y_0)}\{\tau_\infty<t\}=1$ for all $t>0$ if one of the following holds:
(ia)
$r_1\vee0\le \theta_1-1$;
(ib)
$r_2\vee0 \le \theta_2-1$.
\item[{\normalfont(ii)}]
$(X,Y)$ stays infinite if one of the following holds:
(iia)
$r_1\le (\theta_1-1)\vee0$;
(iib)
$r_2\le (\theta_2-1)\vee0$.
\end{itemize}
\etheorem

The following results show that the behavior of explosion and staying infinite for the system $(X,Y)$ given by \eqref{1.1} are affected by the
interaction terms
$a_1X_s^{\theta_1}Y_s^{\kappa_1}$ and $a_2Y_s^{\theta_2}X_s^{\kappa_2}$.

\btheorem\label{tt1.3}
\begin{itemize}
\item[{\normalfont(i)}]
$\lim_{X_0,Y_0\to\infty}\mbf{P}_{(X_0,Y_0)}\{\tau_\infty<t\}=1$ for all $t>0$
if
$\theta_i-1<r_i\vee0$ for $i=1,2$ and
$[(r_1\vee0)+1-\theta_1][(r_2\vee0)+1-\theta_2]<\kappa_1\kappa_2$.
\item[{\normalfont(ii)}] $(X,Y)$ stays infinite if
$(\theta_i-1)\vee0< r_i$  for $i=1,2$ and
 \beqlb\label{1.10}
(r_1+1-\theta_1)(r_2+1-\theta_2)<\kappa_1\kappa_2.
 \eeqlb
\end{itemize}
\etheorem

We state the following condition, which means
the drift terms play a dominant role when the processes are large enough,
and it will be used to consider the critical case.
\bcondition\label{c3}
One of the following holds:
\begin{itemize}
\item[{\normalfont(i)}]
For $b_{10}b_{20}\neq0$,
 \beqlb\label{1.8a}
r_{10}-1>\max_{j=1,2}\{r_{1j}-\varrho_{1j}: b_{1j}\neq0\}
 \eeqlb
and
 \beqlb\label{1.9a}
r_{20}-1>\max_{j=1,2}\{r_{2j}-\varrho_{2j}: b_{2j}\neq0\};
 \eeqlb

\item[{\normalfont(ii)}]
In addition to $\kappa_2 r_2/( r_2+1-\theta_2)< r_1$ for $r_2>0$,
either
\eqref{1.8a} holds and \eqref{1.9a}
is not satisfied for $b_{10}b_{20}\neq0$, or
$b_{20}=0$ and
\eqref{1.8a} holds for $b_{10}\neq0$;
\item[{\normalfont(iii)}]
In addition to $\kappa_1 r_1/( r_1+1-\theta_1)< r_2$ for $r_1>0$,
either
\eqref{1.9a} holds and \eqref{1.8a}
is not satisfied for $b_{10}b_{20}\neq0$, or
$b_{10}=0$ and
\eqref{1.9a} holds for $b_{20}\neq0$.
\end{itemize}
\econdition

\btheorem\label{t1.4}
Suppose that $(\theta_i-1)\vee0< r_i$ for $i=1,2$ and Condition \ref{c3}(i) is satisfied and that
$(r_1+1-\theta_1)(r_2+1-\theta_2)=\kappa_1\kappa_2$.  Then
$\lim_{X_0,Y_0\to\infty}\mbf{P}_{(X_0,Y_0)}\{\tau_\infty<t\}=1$ for all $t>0$,
and $(X,Y)$ stays infinite if
$(a_1/b_1)^{1/( r_1+1-\theta_1)}(a_2/b_2)^{1/\kappa_2}>1$.
\etheorem

Under condition $(r_1+1-\theta_1)(r_2+1-\theta_2)=\kappa_1\kappa_2$, the inequality $(a_1/b_1)^{\frac{1}{r_1+1-\theta_1}}(a_2/b_2)^{\frac{1}{\kappa_2}}>1$ holds if and only if
$(a_2/b_2)^{\frac{1}{r_2+1-\theta_2}}(a_1/b_1)^{\frac{1}{\kappa_1}}>1$,
revealing a symmetric structure.
The following two theorems state the conditions on the nonexplosion and
coming down from infinity
of system \eqref{1.1}.
\btheorem\label{t1.8}
$\mbf{P}\{\tau_\infty<\infty\}=0$
if $\theta_i-1<r_i\vee0$ for $i=1,2$ and one of the following holds:
\begin{itemize}
\item[{\normalfont(i)}]
$r_1,r_2\le0$ and $(1-\theta_1)(1-\theta_2)\ge\kappa_1\kappa_2$;
\item[{\normalfont(ii)}]
$r_1\le0$, $r_2>0$, $\theta_2>0$, $r_2+1>\kappa_1$ and $(1-\theta_1)(r_2+1-\theta_2)\ge\kappa_1\kappa_2$;
\item[{\normalfont(iii)}]
$r_1>0$, $r_2\le0$, $\theta_1>0$, $r_1+1>\kappa_2$ and $(r_1+1-\theta_1)(1-\theta_2)\ge\kappa_1\kappa_2$;
\item[{\normalfont(iv)}]
$r_1,r_2>0$ and at least one of the following holds:
\begin{itemize}
\item[{\normalfont(iva)}]
Condition \ref{c3} is satisfied,
 \beqlb\label{1.7a}
\frac{r_{12}-\alpha_1}{r_2}<\frac{r_1+1-\theta_1}{\kappa_1}\mbox{ if } b_{12}>0,
\qquad \frac{r_{22}-\alpha_2}{r_1}<\frac{r_2+1-\theta_2}{\kappa_2}
\mbox{ if } b_{22}>0
 \eeqlb
and
 \beqlb\label{1.7}
( r_1+1-\theta_1)( r_2+1-\theta_2)>\kappa_1\kappa_2;
  \eeqlb
\item[{\normalfont(ivb)}]
Condition \ref{c3} is satisfied, \eqref{1.7a} holds and
$$(r_1+1-\theta_1)(r_2+1-\theta_2)=\kappa_1\kappa_2,~~(a_1/b_1)^{1/( r_1+1-\theta_1)}(a_2/b_2)^{1/\kappa_2}<1;
 $$
\item[{\normalfont(ivc)}]
\eqref{1.7} holds,
 \beqlb\label{1.6b}
r_{i0}-1\le\max_{j=1,2}\{r_{ij}-\varrho_{ij}:b_{ij}\neq0\}~
\mbox{when}~b_{i0}\neq0
 \eeqlb
for $i=1,2$,
and either
 \beqlb\label{1.5}
(r_1+1-\theta_1)/\kappa_1>[(r_1+1)/(r_2+1)]\vee[r_1/r_2]
 \eeqlb
or
 \beqlb\label{1.6}
(r_2+1-\theta_2)/\kappa_2>[(r_2+1)/(r_1+1)]\vee[r_2/r_1].
 \eeqlb
\end{itemize}
\end{itemize}
\etheorem

\btheorem\label{t1.9}
$(X,Y)$ comes down from infinity if
$\theta_i-1<r_i$ for $i=1,2$
and condition (iv) in Theorem \ref{t1.8} holds.
\etheorem

The conditions in Theorem \ref{t1.8}(iva)--(ivb) consider the case of drift terms dominance, the condition in
Theorem \ref{t1.8}(ivc) considers the case of diffusive terms dominance.
When $r_1= r_2>0$, inequality \eqref{1.7} implies
\eqref{1.5} or \eqref{1.6}.
In what follows, we consider the critical case \eqref{6.2}.
For $i=1,2$ let $\bar{r}_i:=\max_{j=1,2}\{r_{ij}-\varrho_{ij}:b_{ij}>0\}$.
\btheorem\label{t1.111}
Suppose that $(\theta_i-1)\vee0< r_i$ for $i=1,2$ and Condition \ref{c3}(i) is satisfied.
Assume further that
 \beqlb\label{6.2}
(r_1+1-\theta_1)(r_2+1-\theta_2)=\kappa_1\kappa_2,\quad (a_1/b_1)^{1/( r_1+1-\theta_1)}(a_2/b_2)^{1/\kappa_2}=1.
 \eeqlb
Then
\begin{itemize}
\item[{\normalfont(i)}]
$(X,Y)$ stays infinite when
$\bar{r}_i\le0$ for $i=1,2$.
\item[{\normalfont(ii)}]
$(X,Y)$ comes down from infinity when
$\bar{r}_i>0$ for $i=1,2$ and
\beqlb\label{1.11}
 \begin{aligned}
&\frac{r_1+1-\theta_1}{r_2+\theta_2} \le\frac{r_1+1-\theta_1}{\kappa_1}=\frac{\kappa_2}{ r_2+1-\theta_2}
\le\frac{r_1+\theta_2}{r_2+1-\theta_2},\mbox{ or } \\
& \frac{r_1+\theta_2}{r_2+1-\theta_2} \le\frac{r_1+1-\theta_1}{\kappa_1}=\frac{\kappa_2}{ r_2+1-\theta_2}
\le\frac{r_1+1-\theta_1}{r_2+\theta_2}.
 \end{aligned}
 \eeqlb
\item[{\normalfont(iii)}]
$\mbf{P}\{\tau_\infty<\infty\}=0$
either $\bar{r}_i<0$ for $i=1,2$, or
\eqref{1.11} holds.
\end{itemize}
\etheorem

Combining Theorems \ref{tt1.3a}--\ref{tt1.3}, \ref{t1.4} and
\ref{t1.8}--\ref{t1.111}, all eventualities have been accounted for within a specific scenario.

\bremark\label{r1.2}
Let
$ r_1= r_2:=r$, $\theta_1=\theta_2=:\theta$,
$\kappa_1=\kappa_2=:\kappa$ and $\bar{r}_1=\bar{r}_2:=\bar{r}$.
\begin{itemize}
\item[{\normalfont(i)}]
Let $r\le0$.
Then $(X,Y)$ stays infinite.
Moreover,
$\lim_{X_0,Y_0\to\infty}\mbf{P}_{(X_0,Y_0)}\{\tau_\infty<t\}=1$ for all $t>0$ if
$1-\theta<\kappa$,
and $\mbf{P}\{\tau_\infty<\infty\}=0$ if $1-\theta\ge\kappa$.
\item[{\normalfont(ii)}]
Let $r>0$.
If $r<\kappa+\theta-1$, then $\lim_{X_0,Y_0\to\infty}\mbf{P}_{(X_0,Y_0)}\{\tau_\infty<t\}=1$ for all $t>0$
and stays infinite.
If $r>\kappa+\theta-1$,
then
$\mbf{P}\{\tau_\infty<\infty\}=0$
and
$(X,Y)$ comes down from infinity.

\item[{\normalfont(iii)}]
Suppose that $0<r=\kappa+\theta-1$
and Condition \ref{c3}(i) holds. Then
\begin{itemize}
\item[{\normalfont(iiia)}]
$\lim_{X_0,Y_0\to\infty}\mbf{P}_{(X_0,Y_0)}\{\tau_\infty<t\}=1$ for all $t>0$  and $(X,Y)$ stays infinite when  $(a_1/b_1)^{\frac{1}{r_1+1-\theta_1} }(a_2/b_2)^{\frac{1}{\kappa_2}}>1$;
\item[{\normalfont(iiib)}]
$\mbf{P}\{\tau_\infty<\infty\}=0$ and
$(X,Y)$ comes down from infinity if $(a_1/b_1)^{\frac{1}{r_1+1-\theta_1}}(a_2/b_2)^{\frac{1}{\kappa_2}}<1$;
\item[{\normalfont(iiic)}]
under the condition $(a_1/b_1)^{1/( r_1+1-\theta_1)}(a_2/b_2)^{1/\kappa_2}=1$,
we have $\mbf{P}\{\tau_\infty<\infty\}=0$, and
\begin{itemize}
\item
$(X,Y)$ stays infinite if $\bar{r}\le0$
\item
$(X,Y)$ comes down from infinity if $\bar{r}>0$.
\end{itemize}
\end{itemize}
\end{itemize}
\eremark

Before concluding this section, let us briefly outline the proof methods for the main results.
We establish new criteria in Propositions \ref{t2.2} and \ref{t3.1b} for explosion and staying infinite, which are then used to prove Theorems \ref{tt1.3a}--\ref{tt1.3} and \ref{t1.4} and Theorem \ref{t1.111}(i), with certain two-dimensional power functions serving as key test functions (see \eqref{4.5}, \eqref{4.5b}, \eqref{4.5ba}, \eqref{4.8} and \eqref{4.8b}).
We also establish new criteria in Propositions \ref{t3.1} and \ref{t2.4b} for nonexplosion and coming down from infinity, which are employed to prove Theorems \ref{t1.8} and \ref{t1.9} and Theorem \ref{t1.111}(ii)--(iii), where the test functions
$g$ are selected as in Lemmas \ref{tl3.1}--\ref{t2.3} and \eqref{3.6b}.
When the drift term plays a dominant role (Condition \ref{c3}) for sufficiently large system states, the key functions $g$ are given in Lemmas \ref{t8.2} and \ref{t8.3}. For the case where diffusive terms dominate (inequality \eqref{1.6b}), we replace
$y^{\rho_2}$ by $y(1+y)^{-\varepsilon}$
for small enough $\varepsilon>0$ when \eqref{1.5} holds, and we similarly replace
$x^{\rho_1}$ by $x(1+x)^{-\varepsilon}$ under \eqref{1.6} (see Lemma \ref{t8.5}).

Analogous results can be obtained when $b_{10}b_{20}<0$ or $b_{10},b_{20}\le0$;
these are omitted in this paper.
Let $C^2((0,\infty))$ and $C^2((0,\infty)\times(0,\infty))$
denote the second-order continuous differentiable functions spaces on $(0,\infty)$
and $(0,\infty)\times(0,\infty)$, respectively.
The remainder of this paper is organized as follows. In Section \ref{Criteria}, we present criteria for explosion, nonexplosion,
staying infinite and coming down from infinity in two-dimensional processes. The proofs of the theorems are given in Section 3.

\section{Criteria for behaviors near infinity}\label{Criteria}
\setcounter{equation}{0}

In this section we establish some criteria
which will be used to prove the theorems in this paper.
In the following let $((x_t,y_t))_{t\ge0}$ be a two-dimensional
process where $(x_t)_{t\ge0}$ and $(y_t)_{t\ge0}$ are two nonnegative
c\'adl\'ag processes
without negative jumps defined before the minimum of their first times of hitting
zero or explosion.
In this section we always assume that $x_0$ and $y_0$ are deterministic
and let $\mathcal{L}$ denote the operator so that for each $g\in C^2((0,\infty)\times(0,\infty))$  with
 \beqlb\label{3.1a}
\mathcal{L}g \quad
\mbox{ bounded on }
[u,v]\times[u,v] \mbox{~~for all~} v\ge u>0,
 \eeqlb
the process
 \beqlb\label{3.1}
t\mapsto M^g_{t\wedge\gamma_{m,n}} \quad \mbox{is a martingale},
 \eeqlb
where
 \beqnn
M^g_t:=g(x_t,y_t)-g(x_0,y_0)-\int_0^t\mathcal{L}g(x_s,y_s)\dd s
 \eeqnn
and $\gamma_{m,n}:=\tau_m^-\wedge\tau_n^+$
with $\tau_m^-:=\tau^-_m(x)\wedge\tau^-_m(y)$
and $\tau^+_n:=\tau^+_n(x)\wedge\tau^+_n(y)$.
The operator $\mathcal{L}$ can be obtained by It\^o's formula for SDE system.
In this section let $\tau_0:=\tau_0(x)\wedge\tau_0(y)$
and $\tau_\infty:=\lim_{n\to\infty}\tau_n^{+}$.
The following criteria for two-dimensional process
$((x_t,y_t))_{t\ge0}$ generalize Chen's criteria for the uniqueness problem of Markov jump processes.

The following Proposition \ref{t3.1} is the criterion on nonexplosion
and used to prove Theorem \ref{t1.8} and Theorem \ref{t1.111}(iii), and is proved by using \eqref{3.1}
and Gronwall's inequality.
The test function $g$ in using this proposition is often constructed by the
composition of two-dimensional power functions and logarithmic functions.

\bproposition\label{t3.1}
If for each $v>0$, there is a nonnegative function $g_v\in C^2((0,\infty)\times(0,\infty))$ satisfying \eqref{3.1a}
and a constant $d_v>0$ so that
\begin{itemize}
\item[{\normalfont(i)}]
$\lim_{x\vee y\to\infty}g_v(x,y)=\infty$;
\item[{\normalfont(ii)}]
$\mathcal{L}g_v(x,y)\le d_vg_v(x,y)$ for all $x,y\ge v$,
\end{itemize}
then $\mbf{P}\{\tau_\infty<\infty\}=0$.
\eproposition
\proof
Let $v>0$.
By \eqref{3.1} and condition (ii), for all $n\ge1$,
 \beqnn
\mbf{E}\big[g_v(x_{t\wedge\gamma_{v,n}},y_{t\wedge\gamma_{v,n}})\big]
 \ar=\ar
g_v(x_0,y_0)+\mbf{E}\Big[\int_0^{t\wedge\gamma_{v,n}}
\mathcal{L}g_v(x_s,y_s)\dd s\Big] \cr
 \ar\le\ar
g_v(x_0,y_0)+d_v\mbf{E}\Big[\int_0^t
g_v(x_{s\wedge\gamma_{v,n}},y_{s\wedge\gamma_{v,n}})\dd s\Big].
 \eeqnn
It follows from Gronwall's inequality that
 \beqlb\label{3.8}
\mbf{E}\big[g_v(x_{t\wedge\gamma_{v,n}},y_{t\wedge\gamma_{v,n}})\big]
\le g_v(x_0,y_0)\e^{d_vt}.
 \eeqlb
Using Fatou's lemma we obtain
 \beqnn
\mbf{E}\big[g_v(x_{t\wedge \tau_v^-\wedge\tau_\infty},y_{t\wedge  \tau_v^-\wedge \tau_\infty})\big]
 \ar=\ar
\mbf{E}\big[\lim_{n\to\infty}
g_v(x_{t\wedge\gamma_{v,n}},y_{t\wedge\gamma_{v,n}})\big] \cr
 \ar\le\ar
\liminf_{n\to\infty}
\mbf{E}\big[g_v(x_{t\wedge \gamma_{v,n}},y_{t\wedge\gamma_{v,n}})\big]\le g_v(x_0,y_0)\e^{d_vt},
 \eeqnn
which deduces from condition (i) that
$\mbf{P}\{\tau_\infty> t\wedge \tau_v^-\}=1$ for all $t>0$.
Taking $t\to\infty$ one gets
$\mbf{P}\{\tau_\infty\ge \tau_v^-\}=1$,
which gives $\tau_0\le \tau_\infty$ almost surely by taking $v\to0$.
Let $\tau_\infty(x)=\lim_{n\to\infty}\tau^+_n(x)$ and
$\tau_\infty(y)=\lim_{n\to\infty}\tau^+_n(y)$.
Since the processes $(x_t)_{t\ge0}$ and $(y_t)_{t\ge0}$ are defined  before the minimum of their first times of hitting zero or explosion,
we conclude the assertion if $\mbf{P}\{\tau_0=\tau_\infty<\infty\}=0$.
We first assume that $\tau_\infty(x)=\tau_0(y)<\infty$.
Let $h_v(x,y)=g_v(x,v)$ for all $x,y>0$.
By the same argument of above with test function $g$ replaced by $h$, we get
 \beqlb\label{a3.8a}
\mbf{P}\{\tau_\infty(x)> t\wedge\tau_0(y)\wedge\tau_\infty(y) \wedge\tau_v^-(x)\}=1,
\quad t,v>0.
 \eeqlb
Since the two processes are defined  before the minimum of their first times of hitting zero or explosion,
$\tau_0(x)=\tau_\infty(y)=\infty$.
It thus follows from \eqref{a3.8a} that $\tau_0(y)<\tau_\infty(x)$,
which is a contraction.
By the same arguments we have
$\mbf{P}\{\tau_0(x)=\tau_\infty(y)<\infty\}=\mbf{P}\{\tau_0(y)=\tau_\infty(x)<\infty\}=0$,
and then $\mbf{P}\{\tau_0=\tau_\infty<\infty\}$.
This completes the proof.
\qed

The following assertion is the criterion on staying infinite
and used to prove Theorems \ref{tt1.3a}(ii), \ref{tt1.3}(ii), \ref{t1.4} and \ref{t1.111}(i).
The key test function $g$ is often chosen by two-dimensional power functions.
\bproposition\label{t3.1b}
Process $(x,y)$ stays infinite if there is a nonnegative function $g\in C^2((0,\infty)\times(0,\infty))$
satisfying \eqref{3.1a}, $\lim_{u,v\to\infty}g(u,v)=0$ and
 \beqlb\label{3.11}
0<\inf_{v\ge a}g(a,v)\le \sup_{v\ge a}g(a,v)<\infty,\quad
0<\inf_{v\ge a}g(v,a)\le \sup_{v\ge a}g(v,a)<\infty
 \eeqlb
for all $a>0$,
and there is a constant $d_a>0$ so that $\mathcal{L}g(u,v)\le d_ag(u,v)$ for all $u,v\ge a$.
\eproposition
\proof
By \eqref{3.8} we obtain $\mbf{E}[g(x_{t\wedge\gamma_{a,b}},y_{t\wedge\gamma_{a,b}})]
\le g(x_0,y_0)\e^{d_at}$.
It follows from dominated convergence and \eqref{3.11} that
 \beqnn
\mbf{E}\big[g(x_{\tau^-_a},y_{\tau^-_a})1_{\{\tau^-_a<t\wedge\tau_\infty\}}\big]
 \ar\le\ar
\lim_{b\to\infty}\mbf{E}\big[g(x_{\tau^-_a},y_{\tau^-_a})1_{\{\tau^-_a\le t\wedge\tau_b^+\}}\big] \cr
 \ar\le\ar
\lim_{b\to\infty}\mbf{E}\big[g(x_{t\wedge\gamma_{a,b}},y_{t\wedge\gamma_{a,b}})\big]
\le g(x_0,y_0)\e^{d_at}.
 \eeqnn
Since $(x_t)_{t\ge0}$ and $(y_t)_{t\ge0}$ are nonnegative
c\'adl\'ag processes without negative jumps, then
 \beqnn
\big[[\inf_{v\ge a}g(a,v)]\wedge
[\inf_{v\ge a}g(u,a)]\big]\mbf{P}_{(x_0,y_0)}\{\tau^-_a<t\wedge\tau_\infty\}
\le g(x_0,y_0)\e^{d_at}.
 \eeqnn
Using the assumption $\lim_{u,v\to\infty}g(u,v)=0$ we get
$\lim_{x_0,y_0\to\infty}\mbf{P}_{(x_0,y_0)}\{\tau^-_a<t\wedge\tau_\infty\}=0$.
Since the processes $(x_t)_{t\ge0}$ and $(y_t)_{t\ge0}$ are defined  before the minimum of their first times of explosion,
then $\{\tau^-_a\ge \tau_\infty\}\subset\{\tau^-_a=\infty\}$.
Therefore,
 \beqnn
\mbf{P}_{(x_0,y_0)}\{\tau^-_a<t\}
 \ar=\ar
\mbf{P}_{(x_0,y_0)}\{\tau^-_a<t,\tau^-_a< \tau_\infty\}
+
\mbf{P}_{(x_0,y_0)}\{\tau^-_a<t,\tau^-_a\ge \tau_\infty\} \cr
 \ar\le\ar
\mbf{P}_{(x_0,y_0)}\{\tau^-_a<t\wedge\tau_\infty\}.
 \eeqnn
Now letting $x_0,y_0\to\infty$ one concludes the assertion.
\qed

The following Proposition \ref{t2.2} is a criterion
for explosion and used to establish the proofs
of Theorems \ref{tt1.3a}(i) and \ref{tt1.3}(i) and Theorem \ref{t1.4}.
The proof of Proposition \ref{t2.2}, can be shown by \eqref{3.1}
and an martingale technique,
is an adapting the approach for Chen's criteria.
The key test function $g$ in using Proposition \ref{t2.2}
is often selected by two-dimensional power functions.

\bproposition\label{t2.2}
Let $v>0$ and $x_0,y_0>v$ be fixed.
Suppose that there is a function $g\in C^2((0,\infty)\times(0,\infty))$
satisfying \eqref{3.1a} and
bounded on $[v,\infty)\times[v,\infty)$
and a constant $d>0$ such that
\begin{itemize}
\item[{\normalfont(i)}]
$g(x_0,y_0)>0$, and
$g(x,v)\le0$ and $g(v,x)\le0$ for all $x>0$;
\item[{\normalfont(ii)}]
$\mathcal{L}g(x,y)\ge dg(x,y)$ for all $x,y\ge v$.
\end{itemize}
Then $\mbf{P}\{\tau_\infty<\infty\}>0$.
Moreover, $\lim_{x_0,y_0\to\infty}\mbf{P}_{(x_0,y_0)}\{\tau_\infty<t\}=1$
for all $t>0$ if
 \beqnn
\sup_{x,y\ge v}g(x,y)=\lim_{x_0,y_0\to\infty}g(x_0,y_0).
 \eeqnn
\eproposition
\proof
In view of \eqref{3.1}, for all large enough $n\ge1$,
 \beqnn
\mbf{E}\big[g(x_{t\wedge\gamma_{v,n}},y_{t\wedge\gamma_{v,n}})\big]
=
g(x_0,y_0)+\int_0^t\mbf{E}\big[
\mathcal{L}g(x_s,y_s)1_{\{s\le\gamma_{v,n}\}}\big]\dd s.
 \eeqnn
It then follows from integration by parts that
 \beqnn
 \ar\ar
\e^{-dt}\mbf{E}\big[g(x_{t\wedge\gamma_{v,n}},y_{t\wedge\gamma_{v,n}})\big] \cr
 \ar=\ar
g(x_0,y_0)+\int_0^t\e^{-ds}
\dd\Big(\mbf{E}\big[g(x_{s\wedge\gamma_{v,n}},y_{s\wedge\gamma_{v,n}})\big]\Big)
+\int_0^t
\mbf{E}\big[g(x_{s\wedge\gamma_{v,n}},y_{s\wedge\gamma_{v,n}})\big]
\dd\e^{-ds} \cr
 \ar=\ar
g(x_0,y_0)
+\int_0^t\e^{-ds} \mbf{E}\big[
\mathcal{L}g(x_s,y_s)1_{\{s\le\gamma_{v,n}\}}\big] \dd s
-d\int_0^t
\e^{-ds} \mbf{E}\big[g(x_{s\wedge\gamma_{v,n}},y_{s\wedge\gamma_{v,n}})\big]
\dd s.
 \eeqnn
Under condition (ii) we have
 \beqnn
 \ar\ar
\e^{-dt}\mbf{E}\big[g(x_{t\wedge\gamma_{v,n}},y_{t\wedge\gamma_{v,n}})\big]
-g(x_0,y_0) \cr
 \ar\ge\ar
d\int_0^t\e^{-ds} \mbf{E}\Big[
g(x_s,y_s)1_{\{s\le\gamma_{v,n}\}}\big] \dd s
-d\int_0^t
\e^{-ds} \mbf{E}\big[g(x_{s\wedge\gamma_{v,n}},y_{s\wedge\gamma_{v,n}})\big]
\dd s \cr
 \ar=\ar
-d\int_0^t\e^{-ds} \mbf{E}\big[
g(x_{s\wedge\gamma_{v,n}},y_{s\wedge\gamma_{v,n}})1_{\{s>\gamma_{v,n}\}}\big] \dd s.
 \eeqnn
Letting $t\to\infty$ and using the dominated convergence we obtain
 \beqlb\label{3.2}
g(x_0,y_0)
 \ar\le\ar
d\int_0^\infty\e^{-ds} \mbf{E}\big[
g(x_{s\wedge\gamma_{v,n}},y_{s\wedge\gamma_{v,n}})1_{\{s>\gamma_{v,n}\}}\big] \dd s \cr
 \ar=\ar
\mbf{E}\Big[g(x_{\gamma_{v,n}},y_{\gamma_{v,n}})1_{\{\gamma_{v,n}<\infty\}}
d\int_{\gamma_{v,n}}^\infty\e^{-ds}\dd s\Big] \cr
 \ar=\ar
\mbf{E}\big[g(x_{\gamma_{v,n}},y_{\gamma_{v,n}})
\e^{-d\gamma_{v,n}}1_{\{\gamma_{v,n}<\infty\}}\big].
 \eeqlb
Then by the dominated convergence again, we get
 \beqlb\label{3.9}
g(x_0,y_0)
 \ar\le\ar
\lim_{n\to\infty}\mbf{E}\big[g(x_{\gamma_{v,n}},y_{\gamma_{v,n}})
\e^{-d\gamma_{v,n}}1_{\{\gamma_{v,n}<\infty\}}\big] \cr
 \ar=\ar
\mbf{E}\big[g(x_{\tau^-_v\wedge\tau_\infty},y_{\tau^-_v\wedge\tau_\infty})
\e^{-d(\tau^-_v\wedge\tau_\infty)}1_{\{\tau^-_v\wedge\tau_\infty<\infty\}}\big] \cr
 \ar=\ar
\mbf{E}\Big[g(x_{\tau^-_v\wedge\tau_\infty},y_{\tau^-_v\wedge\tau_\infty})
\e^{-d(\tau^-_v\wedge\tau_\infty)}
\big(1_{\{\tau_\infty<\tau^-_v\}}+1_{\{\tau_\infty>\tau^-_v\}}\big)\Big] \cr
 \ar=\ar
 \mbf{E}\Big[g(x_{\tau_\infty},y_{\tau_\infty})
\e^{-d\tau_\infty}
1_{\{\tau_\infty<\tau^-_v\}}
+
g(x_{\tau^-_v},y_{\tau^-_v})
\e^{-d\tau^-_v}1_{\{\tau_\infty>\tau^-_v\}}\Big],
 \eeqlb
which implies
 \beqnn
g(x_0,y_0)\le C_0\mbf{E}\Big[\e^{-d\tau_\infty}
1_{\{\tau_\infty<\tau^-_v\}}\Big]
\le C_0\mbf{P}\{\tau_\infty<\infty\},
 \eeqnn
where $C_0:=\sup_{x,y\ge v}g(x,y)$. This gives the first assertion.
Moreover,
 \beqnn
1=\lim_{x_0,y_0\to\infty}g(x_0,y_0)/ C_0\le \lim_{x_0,y_0\to\infty}\mbf{E}_{(x_0,y_0)}\Big[\e^{-d\tau_\infty}
1_{\{\tau_\infty<\tau^-_v\}}\Big]
 \eeqnn
if
$C_0:=\lim_{x_0,y_0\to\infty}g(x_0,y_0)$.
This implies the second assertion.
\qed

The following proposition gives the criterion
of coming down from infinity and is used to prove Theorem \ref{t1.9}
and Theorem \ref{t1.111}(ii).
\bproposition\label{t2.4b}
Let $z=h(u,v)\in C^2((0,\infty)\times(0,\infty))$ be a nonnegative function such that
\eqref{3.1a} holds,
$u\mapsto h(u,v)$ and $v\mapsto h(u,v)$ are nondecreasing
for all $u,v>0$.
We also assume that there are constants $c_0,\rho_1,\rho_2>0$ such that
$h(u,v)\ge c_0(u^{\rho_1}+v^{\rho_2})$ for all $u,v>a$ and all large enough $a>0$.
Suppose that there are a nonnegative decreasing function
$g\in C^2((0,\infty))$ and function $d: (0,\infty)\mapsto (0,\infty)$ such that
$\lim_{a\to\infty}d(a)=\infty$,
$\lim_{u,v\to\infty}g(h(u,v))=1$, and
 \beqlb\label{3.3b}
\mathcal{L}g(h(u,v))\ge d(a) g(h(u,v)),\qquad h(u,v)\ge a.
 \eeqlb
Then
 \beqnn
\lim_{a\to\infty}\lim_{x_0,y_0\to\infty}\mbf{P}_{(x_0,y_0)}\{\tau_a^-(X)\vee \tau_a^-(Y)<t\}=1,
 \eeqnn
where the process
$z=(z_t)_{t\ge0}$ is defined by $z_t:=h(x_t,y_t)$.
\eproposition
\proof
Let $\bar{\gamma}_{m,n}:=\gamma_{m,n}\wedge \tau_a^-(z)$
and $\bar{\gamma}_{0,\infty}:=\tau_0\wedge\tau_\infty\wedge \tau_a^-(z)$.
Similar to \eqref{3.2} we get
 \beqnn
g(h(x_0,y_0))
\le
\mbf{E}[g(h(x_{\bar{\gamma}_{m,n}},y_{\bar{\gamma}_{m,n}}))
\e^{-d(a)\bar{\gamma}_{m,n}}1_{\{\bar{\gamma}_{m,n}<\infty\}}].
 \eeqnn
Then by dominated convergence,
 \beqnn
 \ar\ar
g(h(x_0,y_0)) \cr
 \ar\le\ar
\mbf{E}\Big[\lim_{m\to0,n\to\infty}g(h(x_{\bar{\gamma}_{m,n}},y_{\bar{\gamma}_{m,n}}))
\e^{-d(a)(\tau_{a}^-(z)\wedge\tau_0\wedge\tau_\infty)}
\big[1_{\{\tau_a^-(z)\le\tau_0\wedge\tau_\infty \}}
+1_{\{\tau_a^-(z)>\tau_0\wedge\tau_\infty \}}\big]\Big] \cr
 \ar\le\ar
\mbf{E}\Big[\lim_{m\to0,n\to\infty}g(h(x_{\bar{\gamma}_{m,n}},y_{\bar{\gamma}_{m,n}}))
\e^{-d(a)(\tau_{a}^-(z)\wedge\tau_0\wedge\tau_\infty)}
\big[1_{\{\tau_a^-(z)\le\tau_0\wedge\tau_\infty ,\tau_a^-(z)<t\}} \cr
 \ar\ar\qquad
+1_{\{\tau_a^-(z)\le\tau_0\wedge\tau_\infty ,\tau_a^-\ge t\}}
+1_{\{\tau_a^-(z)>\tau_0\wedge\tau_\infty \}}\big]\Big] \cr
 \ar\le\ar
g(a)\mbf{P}\{\tau_a^-(z)<t\}+g(a)\e^{-d(a)t}+
\mbf{P}\{\tau_a^-(z)>\tau_0\wedge\tau_\infty\}.
 \eeqnn
Since the processes $((x_t,y_t))_{t\ge0}$ is defined before $\tau_0\wedge\tau_\infty$, we have
$\{\tau_a^-(z)>\tau_0\wedge\tau_\infty\}=\emptyset$.
Thus
 \beqnn
g(h(x_0,y_0))
\le
g(a)\mbf{P}\{\tau_a^-(z)<t\}+g(a)\e^{-d(a)t}.
 \eeqnn
Taking $x_0,y_0\to\infty$ and then $a\to\infty$ we get
 \beqnn
\lim_{a\to\infty}\lim_{x_0,y_0\to\infty}\mbf{P}_{(x_0,y_0)}\{\tau_a^-(z)<t\}= 1.
 \eeqnn
Let $a>0$ be large enough.
Since $h(u,v)\ge c_0(u^{\rho_1}+v^{\rho_2})$ for $u,v>a$,
$\{\tau_b^-(X)\vee \tau_b^-(Y)\le\tau_a^-(z)\}$
for $b=(a/c_0)^{\rho_1^{-1}}\vee (a/c_0)^{\rho_2^{-1}}$.
Thus one concludes the assertion.
\qed

\section{Proofs of the main results}
\setcounter{equation}{0}

In this section we state the proofs for assertions in Section 1.
We first state some notations and two lemmas which will be used in the
rest of this section.
For any $g\in C^2((0,\infty)\times(0,\infty))$ and $x,y,z\ge0$ define
 \beqlb\label{3.21}
K_z^1g(x,y):=
g(x+z,y)-g(x,y)-g'_x(x,y)z
 \eeqlb
and
 \beqlb\label{3.22}
K_z^2g(x,y):=
g(x,y+z)-g(x,y)-g'_y(x,y)z.
 \eeqlb
By It\^o's formula, it is easy to see that the operator $\mathcal{L}$
is given by
 \beqlb\label{3.0}
\mathcal{L}g(x,y)
 \ar:=\ar
a_1x^{\theta_1}y^{\kappa_1}g'_x(x,y)
+a_2y^{\theta_2}x^{\kappa_2}g'_y(x,y) \cr
 \ar\ar
-b_{10}x^{r_{10}}g'_x(x,y)
+ b_{11}x^{r_{11}}g''_{xx}(x,y)
+b_{12}x^{r_{12}}\int_0^\infty K_z^1g(x,y)\mu_1(\dd z) \cr
 \ar\ar
-b_{20}y^{r_{20}}g'_y(x,y)
+ b_{21}y^{r_{21}}g''_{yy}(x,y)
+b_{22}y^{r_{22}}\int_0^\infty K_z^2g(x,y)\mu_2(\dd z).
 \eeqlb
For $i=1,2$ and $\rho<\alpha_i$ let $c_i(\rho):=\frac{\Gamma(\alpha_i-\rho)}{
\Gamma(\alpha_i)\Gamma(2-\rho)}$.
By \cite[Lemma 11]{XYZh24} and Taylor's formula,
for $\rho\in(-\infty,0)\cup(0,1)\cup(1,\alpha_i)$ and $i=1,2$,
 \beqlb\label{3.0a}
c_i(\rho):=[\rho(\rho-1)]^{-1}\int_0^\infty [(1+z)^\rho-1-\rho z]\mu_i(\dd z)
=\int_0^\infty z^2\mu_i(\dd z)\int_0^1(1+zu)^{\rho-2}(1-u)\dd u.
 \eeqlb
In the following we state two lemmas which can be gotten from
\cite{XYZh25} and will be used in the proofs.

\blemma\label{A1.1}
\begin{itemize}
\item[{\normalfont(i)}]
Suppose that $p_1,p_2,p_3,p_4>0$ and $c_1,c_2,c_3>0$.
If $p_3/p_1+p_4/p_2<1$,
then there is a constant $c>1$ so that
$c_1x^{p_1}+c_2y^{p_2}
\ge c_3x^{p_3}y^{p_4}$ for all $x,y>c$.
\item[{\normalfont(ii)}]
For $x,y\ge0$, we have
$x^p+y^p\ge (x+y)^p$ for $0<p\le 1$ and
$x^p+y^p\ge 2^{1-p}(x+y)^p$ for $p>1$.
\item[{\normalfont(iii)}]
Let $p,q>1$ satisfy $1/p+1/q\le1$.
Then for each $\delta>0$, there is a constant $C_\delta>0$ so that
$u+v\ge C_\delta
u^{1/p}v^{1/q}$ for any $u,v\ge \delta$.
\end{itemize}
\elemma
\proof
The assertion (ii) is given by \cite[Lemma 3.1(ii)]{XYZh25}
and the assertion (i) is modification of that of \cite[Lemma 3.1(iii)]{XYZh25}.
We obtain the assertion (iii) by \cite[Lemma 3.1(i)]{XYZh25}.
\qed

\blemma\label{t8.1}
\begin{itemize}
\item[{\normalfont(i)}]
If $r_i>\theta_i-1$ for $i=1,2$ and the constants $\rho_1,\rho_2>1$ satisfy
 \beqlb\label{8.1}
\frac{r_1+1-\theta_1}{\kappa_1}
>
\frac{r_1+\rho_1}{r_2+\rho_2}
>\frac{\kappa_2}{r_2+1-\theta_2},
 \eeqlb
then for any constants $c_1,c_2,c_3,c_4>0$ there is a constant $c>1$ so that for all $x,y\ge c$,
 \beqnn
 \ar\ar
c_1x^{r_1+\rho_1}+c_2y^{r_2+\rho_2}
-c_3x^{\theta_1-1+\rho_1}y^{\kappa_1}
-c_4y^{\theta_2-1+\rho_2}x^{\kappa_2}\ge0.
 \eeqnn
\item[{\normalfont(ii)}]
If the constants $\rho_1,\rho_2>\theta_1\vee\theta_2$ satisfy
 \beqnn
\frac{ r_1+1-\theta_1}{\kappa_1}
<
\frac{1+\rho_1+\kappa_2-\theta_1}{1+\rho_2+\kappa_1-\theta_2},
 \eeqnn
then for each $\delta>0$, there
is a constant $c_\delta>1$ such that
 \beqnn
x^{\theta_1-1-\rho_1}y^{\kappa_1}
+y^{\theta_2-1-\rho_2}x^{\kappa_2}
\ge \delta x^{ r_1-\rho_1},\qquad x,y\ge c_\delta.
 \eeqnn
\item[{\normalfont(iii)}]
Given $0<\rho<1$, $\rho_1>0$, let
 \beqlb\label{3.7bb}
K(v,z):=-\big[v[(1+z)^{\rho_1}-1]+1\big]^\rho+1+zv\rho\rho_1,\quad 0\le v\le 1, \, z>0.
 \eeqlb
Then
\begin{itemize}
\item[{\normalfont(iiia)}]
for all $\rho_1\ge2$
and $0<\rho<1/\rho_1$ we have
$\int_0^\infty K(v,z)\mu_1(\dd z)
\ge- \rho \tilde{d} v^{\alpha_1/\rho_1}$
for all $0\le v\le 1$,
where $\tilde{d}=\rho_1^{\alpha_1}\int_0^\infty (u\wedge u^2)\mu_1(\dd u)$.
\item[{\normalfont(iiib)}]
for all $\rho_1>1$, $0\le v\le1$ and $0<\rho<1/\rho_1$, and all $\delta>0$,
 \beqnn
\int_0^\infty K(v,z)\mu_1(\dd z)\ge
\rho\rho_1(1-\rho\rho_1) c_1(\rho\rho_1)v^2
-\rho\rho_1(\rho_1-1)[d_{2,\delta}v^\rho+v(1-v)d_{1,\delta}]
 \eeqnn
for functions
$d_{1,\delta},d_{2,\delta}>0$ satisfying $\lim_{\delta\to\infty}d_{2,\delta}=0$;
\item[{\normalfont(iiic)}]
for all $0<\rho,\rho_1<1$, we have
$\int_0^\infty K(v,z)\mu_1(\dd z)\ge0$ for $0\le v\le 1$.
\end{itemize}
\end{itemize}
\elemma
\proof
The conclusion (i) is obtained by Lemma \ref{A1.1}(i).
The proof of (ii) is obvious when $r_1\le \theta_1-1$,
and is modification to that of \cite[Lemma 3.2]{XYZh25}
when $r_1>\theta_1-1$.
The assertions in (iiib) and (iiic) follow from \cite[Lemma 3.9(ii)]{XYZh25}.
We show the assertion (iiia) in the following.
Observe that
 \beqnn
K_{z}'(v,z)=\rho\rho_1 v[1-(1+z)^{\rho_1-1}H(v,z)^{\rho-1}]
\ge-\rho\rho_1 v (1+z)^{\rho_1-1}H(v,z)^{\rho-1}
 \eeqnn
and
 \beqnn
K_{zz}''(v,z)\ge
-\rho\rho_1(\rho_1-1)
v(1+z)^{\rho_1-2}H(v,z)^{\rho-1},
 \eeqnn
where $H(v,z):=v[(1+z)^{\rho_1}-1]+1\ge 1$.
For $i=1,2$ and $0\le v\le 1$, we have
 \beqnn
v(1+z)^{\rho_1-i}H(v,z)^{\rho-1}
=
v^{i/\rho_1}[v(1+z)^{\rho_1}]^{1-i/\rho_1}H(v,z)^{\rho-1}
\le
v^{i/\rho_1}H(v,z)^{\rho-i/\rho_1}
\le v^{i/\rho_1}.
 \eeqnn
Thus, for $0\le v\le 1$,
 \beqnn
K_{z}'(v,z)\ge- \rho\rho_1v^{1/\rho_1},\qquad
K_{zz}''(v,z)\ge-\rho\rho_1^2v^{2/\rho_1}.
 \eeqnn
It then follows from Taylor's formula that for all $0\le v\le 1$,
 \beqnn
K(v,z)\ar=\ar
z\int_0^1K_{z}'(v,zu)\dd u
\ge- z\rho\rho_1v^{1/\rho_1}, \cr
K(v,z)\ar=\ar z^2\int_0^1K''_{zz}(v,uz)(1-u)\dd u
\ge
-z^2\rho\rho_1^2v^{2/\rho_1}.
 \eeqnn
With the variable $\rho_1v^{1/\rho_1}z$ replaced by $u$ we obtain
 \beqnn
-\int_0^\infty K(v,z)\mu_1(\dd z)
\le \rho\int_0^\infty [(\rho_1v^{1/\rho_1}z)
\wedge(\rho_1v^{1/\rho_1}z)^2]\mu_1(\dd z)
=\rho\rho_1^{\alpha_1}v^{\alpha_1/\rho_1}\int_0^\infty (u\wedge u^2)\mu_1(\dd u),
 \eeqnn
which completes the proof. \qed

\subsection{Proofs of Theorems \ref{tt1.3a}--\ref{tt1.3} and \ref{t1.4}.}

Before proving
Theorems \ref{tt1.3a}--\ref{tt1.3} and \ref{t1.4}, we introduce the following lemma
whose proof is given by Lemmas \ref{A1.1}(ii) and \ref{t8.1}(ii).
The following Lemma \ref{t5.4}(iii) is the key to show Theorem \ref{t1.111}(i).

\blemma\label{t5.4}
\begin{itemize}
\item[{\normalfont(i)}]
Given $0<\rho<1$ and $\rho_1,\rho_2>0$, define a continuous function
 \beqlb\label{4.5}
h(x,y):=(x^{-\rho_1}+y^{-\rho_2})^\rho, \qquad x,y>0.
 \eeqlb
\begin{itemize}
\item[{\normalfont(ia)}]
If there are constants $0<\rho<1$ and $\rho_1,\rho_2>\theta_1\vee\theta_2\vee(r_1-\kappa_2)\vee(r_2-\kappa_1)$ so that
 \beqlb\label{5.1}
\frac{r_1+1-\theta_1}{\kappa_2+\rho_1+1-\theta_1}
<\frac{\kappa_1}{\kappa_1+\rho_2+1-\theta_2},
~~
\frac{r_2+1-\theta_2}{\kappa_1+\rho_2+1-\theta_2}
<\frac{\kappa_2}{\kappa_2+\rho_1+1-\theta_1},
 \eeqlb
and
 \beqlb\label{5.2}
\frac{\rho\rho_1+1-\theta_1}{\kappa_2+\rho_1+1-\theta_1}
<\frac{\kappa_1}{\kappa_1+\rho_2+1-\theta_2},
\quad
\frac{\rho\rho_2+1-\theta_2}{\kappa_1+\rho_2+1-\theta_2}
<\frac{\kappa_2}{\kappa_2+\rho_1+1-\theta_1},
 \eeqlb
then there are constants $c>0$ and $c_1>1$ so that
$\mathcal{L}h(x,y)\le -c$ for all $x,y\ge c_1$.
\item[{\normalfont(ib)}]
If for each $i=1,2$, either $r_i\le0$ or the term
$r_i$ in inequality \eqref{5.1} holds for some constants $\rho_1,\rho_2>\theta_1\vee\theta_2$,
then, for each $c_2>0$, there is a constant $c_3>0$
such that $\mathcal{L}h(x,y)\le c_3 h(x,y)$ for all $x,y\ge c_2$.
\end{itemize}
\item[{\normalfont(ii)}]
Suppose that $r_1,r_2>0$ and Condition \ref{c3}(i) holds and that there is a constant $\varepsilon\in(0,1/2)$ so that
 \beqlb\label{5.3}
\Big[\frac{a_1}{b_1(1+\varepsilon)}\Big]^{1/(r_1+1-\theta_1)}
\Big[\frac{a_2}{b_2(1+\varepsilon)}\Big]^{1/\kappa_2}\ge1.
 \eeqlb
Let $\rho_1,\rho_2>\theta_1\vee\theta_2\vee(r_1-\kappa_2)\vee(r_2-\kappa_1)$ be the constants satisfy
 \beqlb\label{5.4}
\frac{r_1+1-\theta_1}{\kappa_1}
= \frac{\kappa_2}{r_2+1-\theta_2}
=
\frac{1+\rho_1+\kappa_2-\theta_1}{1+\rho_2+\kappa_1-\theta_2}.
 \eeqlb
Then there are constants $\delta_0,\tilde{c}>0$, $0<\rho<1$ and $\tilde{c}_1>1$ such that
$\mathcal{L}h(x,y)\le -\tilde{c}$ for all $x,y\ge \tilde{c}_1$,
where function $h$ is defined by
 \beqlb\label{4.5b}
h(x,y):=(\delta_0x^{-\rho_1}+y^{-\rho_2})^\rho,\qquad x,y>0.
 \eeqlb
Moreover, for each $c_4>0$, there is a constant $c_5>0$
such that $\mathcal{L}h(x,y)\le c_5 h(x,y)$ for all $x,y\ge c_4$.
\item[{\normalfont(iii)}]
Assume that $r_1,r_2>0$ and Condition \ref{c3}(i) and \eqref{6.2} are fulfilled.
In addition, suppose that \eqref{5.3} holds with $\varepsilon=0$
and that \eqref{5.4} is satisfied for some constants $\rho_1,\rho_2>\theta_1\vee\theta_2\vee(r_1-\kappa_2)\vee(r_2-\kappa_1)$.
If $\bar{r}_i\le0$ for $i=1,2$,
then for each $c_6>0$, there are constants $\delta_0,c_7>0$ such that
$\mathcal{L}h(x,y)\le c_7 h(x,y)$ for all $x,y\ge c_6$,
where
 \beqlb\label{4.5ba}
h(x,y):=\ln (1+\delta_0x^{-\rho_1}+y^{-\rho_2}),\qquad x,y>0.
 \eeqlb
\end{itemize}
\elemma
\proof
(i)
Observe that
 \beqnn
h'_x(x,y)=-\rho\rho_1 (x^{-\rho_1}+y^{-\rho_2})^{\rho-1}x^{-\rho_1-1},~
h''_{xx}(x,y)\le
\rho\rho_1(\rho_1+1) (x^{-\rho_1}+y^{-\rho_2})^{\rho-1}x^{-\rho_1-2}
 \eeqnn
and
 \beqnn
h'_y(x,y)=-\rho\rho_2 (x^{-\rho_1}+y^{-\rho_2})^{\rho-1}y^{-\rho_2-1},~
h''_{yy}(x,y)\le
\rho\rho_2(\rho_2+1) (x^{-\rho_1}+y^{-\rho_2})^{\rho-1}y^{-\rho_2-2}.
 \eeqnn
Since $\mu_1(\dd (vz))=v^{-\alpha_1}\mu_1(\dd z)$,
we first substitute $z$ with $xu$, and then substitute
$u$ with $z$, yielding
 \beqlb\label{4.1a}
 \ar\ar
\int_0^\infty K_z^1h(x,y)\mu_1(\dd z)
=x^{-\alpha_1}\int_0^\infty K_{xu}^1h(x,y)\mu_1(\dd u) \cr
 \ar=\ar
x^{-\alpha_1}(x^{-\rho_1}+y^{-\rho_2})^{\rho}
\int_0^\infty K(x^{-\rho_1}(x^{-\rho_1}+y^{-\rho_2})^{-1},z)\mu_1(\dd z),
 \eeqlb
where
 \beqnn
K(v,z):=\big[v[(1+z)^{-\rho_1}-1]+1\big]^\rho-1+z\rho\rho_1v,\qquad
z\ge0,~0\le v\le1.
 \eeqnn
For $0\le v\le1$,
 \beqnn
K''_{zz}(v,z)\le
v\rho\rho_1(\rho_1+1)\big[v[(1+z)^{-\rho_1}-1]+1\big]^{\rho-1}(1+z)^{-\rho_1-2}
\le v\rho\rho_1(\rho_1+1)(1+z)^{-\rho\rho_1-2}.
 \eeqnn
It thus follows from Taylor's formula and \eqref{3.0a} that
 \beqlb\label{3.10}
 \ar\ar
\int_0^\infty K(v,z)\mu_1(\dd z)
 =
\int_0^\infty z^2\mu_1(\dd z)\int_0^1K''_{zz}(v,zu)(1-u)\dd u \cr
 \ar\le\ar
v\rho\rho_1(\rho_1+1)\int_0^\infty z^2\mu_1(\dd z)
\int_0^1(1+zu)^{-\rho\rho_1-2}(1-u)\dd u
=v\rho\rho_1(\rho_1+1)c_1(-\rho\rho_1),
 \eeqlb
where the function $c_1$ is defined in \eqref{3.0a}.
By \eqref{4.1a},
 \beqnn
\int_0^\infty K_z^1h(x,y)\mu_1(\dd z)
\le\rho\rho_1(\rho_1+1)c_1(-\rho\rho_1)
x^{-\rho_1-\alpha_1}(x^{-\rho_1}+y^{-\rho_2})^{\rho-1}.
 \eeqnn
Similarly,
 \beqnn
\int_0^\infty K_z^2h(x,y)\mu_2(\dd z)
\le\rho\rho_2(\rho_2+1)c_2(-\rho\rho_2)
y^{-\rho_2-\alpha_2}(x^{-\rho_1}+y^{-\rho_2})^{\rho-1}.
 \eeqnn
Now by \eqref{3.0},
there are constants $C_1,C_2>0$ so that
 \beqlb\label{4.2}
 \ar\ar
\rho^{-1}[x^{-\rho_1}+y^{-\rho_2}]^{(1-\rho)}\mathcal{L}h(x,y) \cr
 \ar\le\ar
-a_1\rho_1x^{\theta_1-1-\rho_1}y^{\kappa_1}
-a_2\rho_2y^{\theta_2-1-\rho_2}x^{\kappa_2}
+b_{10}\rho_1x^{-\rho_1+r_{10}-1}  \cr
 \ar\ar
+b_{11} \rho_1(\rho_1+1)x^{-\rho_1+r_{11}-2}
+b_{12} \rho_1(\rho_1+1)c_1(-\rho\rho_1)x^{-\rho_1+r_{12}-\alpha_1} \cr
 \ar\ar
+b_{20}\rho_2y^{-\rho_2+r_{20}-1}
+b_{21} \rho_2(\rho_2+1)y^{-\rho_2+r_{21}-2}
+b_{22} \rho_2(\rho_2+1)c_2(-\rho\rho_2)y^{-\rho_2+r_{22}-\alpha_2} \cr
 \ar\le\ar
[C_1x^{ r_1-\rho_1}
+C_2y^{r_2-\rho_2}]
-[a_1\rho_1x^{\theta_1-1-\rho_1}y^{\kappa_1}
+a_2\rho_2y^{\theta_2-1-\rho_2}x^{\kappa_2}] \cr
 \ar=:\ar
M_1(x,y)-2M_2(x,y),\qquad x,y\ge1.
 \eeqlb
Under \eqref{5.1}, by Lemma \ref{t8.1}(ii),
$M_1(x,y)-M_2(x,y)\le 0$
for all $x,y\ge \tilde{c}_0$ and some constant $\tilde{c}_0>0$.
Moreover, by \eqref{5.2} and Lemmas \ref{t8.1}(ii) and \ref{A1.1}(ii),
there is a constant $c_1>\tilde{c}_0$ so that
 \beqnn
M_2(x,y)\ge x^{-\rho_1(1-\rho)}+y^{-\rho_2(1-\rho)}
\ge [x^{-\rho_1}+y^{-\rho_2}]^{(1-\rho)},\qquad x,y\ge c_1.
 \eeqnn
It thus follows from \eqref{4.2} that
$\mathcal{L}h(x,y)\le-\rho[x^{-\rho_1}+y^{-\rho_2}]^{\rho-1}M_2(x,y)\le -\rho$ for all $x,y\ge c_1$, which concludes
(ia).
When $r_i\le0$, we have $x^{r_i-\rho_i}\le c_2^{r_i}x^{-\rho_i}$ for all $x\ge c_2>0$.
Then by Lemma \ref{t8.1}(ii) again, there is a constant $c_0>1$ so that
for all $x,y\ge c_0$ we have
$M_1(x,y)-M_2(x,y)\le (C_1\vee C_2)[x^{-\rho_1}+y^{-\rho_2}]$
for all $x,y\ge c_0$.
Since $\rho_1,\rho_2>(r_1-\kappa_2)\vee(r_2-\kappa_1)$,
then $M_1(x,y)-M_2(x,y)$  is bounded above on
$([c_2,\infty)\times[c_2,\infty))
\setminus([c_0,\infty)\times[c_0,\infty))$.
Thus there is a constant $C_3>0$ such that
$\mathcal{L}h(x,y)\le C_3h(x,y)$
for all $x,y\ge c_2$,
which implies the assertion (ib) by \eqref{4.2}.

(ii)
Recall functions $c_i$ from \eqref{3.0a}.
Under Condition \ref{c3}(i), by the same arguments as in \eqref{4.2},
 \beqlb\label{4.2b}
 \ar\ar
\rho^{-1}[\delta_0x^{-\rho_1}+y^{-\rho_2}]^{(1-\rho)}\mathcal{L}h(x,y) \cr
 \ar\le\ar
-\delta_0a_1\rho_1x^{\theta_1-1-\rho_1}y^{\kappa_1}
-a_2\rho_2y^{\theta_2-1-\rho_2}x^{\kappa_2}
+\delta_0\rho_1 b_{10}x^{-\rho_1+r_{10}-1}  \cr
 \ar\ar
+\delta_0b_{11} \rho_1(\rho_1+1)x^{-\rho_1+r_{11}-2}
+\delta_0b_{12} \rho_1(\rho_1+1)c_1(-\rho\rho_1)x^{-\rho_1+r_{12}-\alpha_1} \cr
 \ar\ar
+b_{20}\rho_2y^{-\rho_2+r_{20}-1}
+b_{21} \rho_2(\rho_2+1)y^{-\rho_2+r_{21}-2}
+b_{22} \rho_2(\rho_2+1)c_2(-\rho\rho_2)y^{-\rho_2+r_{22}-\alpha_2} \cr
 \ar=\ar
-\tilde{M}_1(x,y)+(1+\varepsilon)\tilde{M}_2(x,y)-\varepsilon \tilde{M}_2(x,y)+\tilde{M}_3(x,y),\qquad x,y>0,
 \eeqlb
where
 \beqnn
\tilde{M}_1(x,y):=
\delta_0a_1\rho_1x^{\theta_1-1-\rho_1}y^{\kappa_1}
+a_2\rho_2y^{\theta_2-1-\rho_2}x^{\kappa_2},~
\tilde{M}_2(x,y):=\delta_0\rho_1 b_1x^{-\rho_1+r_1}
+b_2\rho_2y^{-\rho_2+r_2}
 \eeqnn
and
 \beqnn
\tilde{M}_3(x,y)
 \ar:=\ar
\delta_0b_{11} \rho_1(\rho_1+1)x^{-\rho_1+r_{11}-2}
+\delta_0b_{12} \rho_1(\rho_1+1)c_1(-\rho\rho_1)x^{-\rho_1+r_{12}-\alpha_1} \cr
 \ar\ar
+b_{21} \rho_2(\rho_2+1)y^{-\rho_2+r_{21}-2}
+b_{22} \rho_2(\rho_2+1)c_2(-\rho\rho_2)y^{-\rho_2+r_{22}-\alpha_2}.
 \eeqnn
By \eqref{5.3}--\eqref{5.4} and \cite[Lemma 3.3]{XYZh25},
there is a constant $\delta_0>0$ such that
$-\tilde{M}_1(x,y)+(1+\varepsilon)\tilde{M}_2(x,y)\le0$ for
all $x,y>0$.
Under Condition \ref{c3}(i), there is a constant $\tilde{c}_1>1$
so that
$\tilde{M}_3(x,y)\le2^{-1}\varepsilon \tilde{M}_2(x,y)$
for all $x,y\ge \tilde{c}_1$.
Let $0<\rho<1$ be small enough so that $r_1>\rho\rho_1$ and $r_2>\rho\rho_2$.
Then by Lemma \ref{A1.1}(ii),
 \beqnn
\tilde{M}_2(x,y)
 \ar=\ar
\delta_0\rho_1 b_1x^{r_1-\rho\rho_1} x^{-\rho_1(1-\rho)}
+b_2\rho_2y^{r_2-\rho\rho_2}y^{-\rho_2(1-\rho)}
\ge (b_1\wedge b_2)[\delta_0x^{-\rho_1(1-\rho)}+y^{-\rho_2(1-\rho)}] \cr
 \ar\ge\ar
(b_1\wedge b_2) [\delta_0x^{-\rho_1}+y^{-\rho_2}]^{1-\rho},\qquad x,y\ge1.
 \eeqnn
It thus follows from \eqref{4.2b} that $\mathcal{L}h(x,y)\le -2^{-1}\varepsilon\rho(b_1\wedge b_2)$
for all $x,y\ge \tilde{c}_1$, which states the first assertion.
Using the assumptions $\rho_1,\rho_2>(r_1-\kappa_2)\vee(r_2-\kappa_1)$,
$\tilde{M}_2(x,y)-\tilde{M}_1(x,y)$ has an upper bound on
$([c_4,\infty)\times[c_4,\infty))
\setminus([\tilde{c}_1,\infty)\times[\tilde{c}_1,\infty))$.
Thus there is a there is a constant $C_4>0$ such that
 \beqnn
\tilde{M}_2(x,y)-\tilde{M}_1(x,y)+\tilde{M}_3(x,y)\le C_4[x^{-\rho_1}+y^{-\rho_2}]^{(1-\rho)}
 \eeqnn
for $(x,y)\in([c_4,\infty)\times[c_4,\infty))
\setminus([\tilde{c}_1,\infty)\times[\tilde{c}_1,\infty))$,
which implies the second assertion by \eqref{4.2b} and the first assertion.

(iii)
Let $J(x,y):=1+\delta_0x^{-\rho_1}+y^{-\rho_2}$.
Then
 \beqnn
h'_x(x,y)=-\delta_0\rho_1 J(x,y)^{-1}x^{-\rho_1-1},~
h''_{xx}(x,y)\le
\delta_0\rho_1(\rho_1+1) J(x,y)^{-1}x^{-\rho_1-2}
 \eeqnn
and
 \beqnn
h'_y(x,y)=-\rho_2 J(x,y)^{-1}y^{-\rho_2-1},~
h''_{yy}(x,y)\le
\rho_2(\rho_2+1) J(x,y)^{-1}y^{-\rho_2-2}.
 \eeqnn
As the same arguments in \eqref{4.1a},
 \beqlb\label{4.1}
 \ar\ar
\int_0^\infty K_z^1h(x,y)\mu_1(\dd z)
=x^{-\alpha_1}\int_0^\infty K_{xu}^1h(x,y)\mu_1(\dd u) \cr
 \ar=\ar
x^{-\alpha_1}
\int_0^\infty K(\delta_0x^{-\rho_1}J(x,y)^{-1},z)\mu_1(\dd z),
 \eeqlb
where
 \beqnn
K(v,z):=\ln\big[v[(1+z)^{-\rho_1}-1]+1\big]+z\rho_1v,\qquad
z\ge0,~0\le v\le1.
 \eeqnn
For $0\le v\le1$,
 \beqnn
K''_{zz}(v,z)\le
v\rho_1(\rho_1+1)\big[v[(1+z)^{-\rho_1}-1]+1\big]^{-1}(1+z)^{-\rho_1-2}
\le v\rho_1(\rho_1+1)(1+z)^{-2}.
 \eeqnn
Then using Taylor's formula,
 \beqlb\label{3.10}
 \ar\ar
\int_0^\infty K(v,z)\mu_1(\dd z)
 =
\int_0^\infty z^2\mu_1(\dd z)\int_0^1K''_{zz}(v,zu)(1-u)\dd u \cr
 \ar\le\ar
v\rho_1(\rho_1+1)\int_0^\infty z^2\mu_1(\dd z)
\int_0^1(1+zu)^{-2}(1-u)\dd u
=:v\rho_1(\rho_1+1)\tilde{c}_1.
 \eeqlb
By \eqref{4.1},
 \beqnn
\int_0^\infty K_z^1h(x,y)\mu_1(\dd z)
\le \rho_1(\rho_1+1)\tilde{c}_1\delta_0
x^{-\rho_1-\alpha_1}J(x,y)^{-1}.
 \eeqnn
Similarly, there is a constant $\tilde{c}_2>0$ such that
 \beqnn
\int_0^\infty K_z^2h(x,y)\mu_2(\dd z)
\le \rho_2(\rho_2+1)\tilde{c}_2
y^{-\rho_2-\alpha_2}J(x,y)^{-1}.
 \eeqnn
By \eqref{3.0},
 \beqlb\label{4.2c}
J(x,y)\mathcal{L}h(x,y)
\le
-\tilde{M}_1(x,y)+\tilde{M}_2(x,y)+\tilde{M}_3(x,y),\qquad x,y>0,
 \eeqlb
where the functions $\tilde{M}_1,\tilde{M}_2,\tilde{M}_3$ are defined
in (ii) with $c_i(-\rho\rho_i)$ replaced by
$\tilde{c}_i$ for $i=1,2$. By \cite[Lemma 3.3]{XYZh25},
there is a constant $\delta_0>0$ such that
$\tilde{M}_2(x,y)\le \tilde{M}_1(x,y)$ for
all $x,y>0$.
Since $\bar{r}_i\le0$ for $i=1,2$ under the assumptions,
then for each $c_6>0$, there are constants $C_1,C_2,C_3,C_4>0$ such that
 \beqnn
\tilde{M}_3(x,y)
 \ar\le\ar
C_1 x^{\bar{r}_1-\rho_1}
+C_2 y^{\bar{r}_2-\rho_2}
\le
C_1 c_6^{\bar{r}_1}x^{-\rho_1}
+C_2c_6^{\bar{r}_2} y^{-\rho_2} \cr
 \ar\le\ar
C_3 [\delta_0x^{-\rho_1}
+y^{-\rho_2 }]
\le
C_4 h(x,y),\qquad x,y\ge c_6,
 \eeqnn
which follows from \eqref{4.2c} and the fact $J(x,y)\ge1$
for all $x,y>0$ that for each $c_6>0$, we have
$\mathcal{L}h(x,y)\le C_4 h(x,y)$ for all $x,y\ge c_6$.
This completes the proof.
\qed

Based on Lemma \ref{t5.4}, we can prove
Theorems \ref{tt1.3a}--\ref{tt1.3} and \ref{t1.4}.

\noindent{\it Proof of Theorems \ref{tt1.3a} and \ref{tt1.3}.}
Let $h$ be the function defined in \eqref{4.5} and
 \beqlb\label{4.8}
g(x,y):=[c_1^{-\rho_1}\wedge c_1^{-\rho_2}]^\rho-h(x,y),
 \eeqlb
where $c_1>1$ is the constant determined by Lemma \ref{t5.4}(ia).
Under the assumptions in Theorems \ref{tt1.3a}(i) or \ref{tt1.3}(i),
there are constants  $\rho_1,\rho_2>\theta_1\vee\theta_2$
such that
 \beqnn
(r_1\vee0)+1-\theta_1
<\frac{\kappa_1(\kappa_2+\rho_1+1-\theta_1)}{\kappa_1+\rho_2+1-\theta_2},\quad
(r_2\vee0)+1-\theta_2
<\frac{\kappa_2(\kappa_1+\rho_2+1-\theta_2)}{\kappa_2+\rho_1+1-\theta_1},
 \eeqnn
and then there is a constant $0<\rho<1$ such that
 \eqref{5.1}--\eqref{5.2} holds. Thus one gives conclusions of
Theorems \ref{tt1.3a}(i) and \ref{tt1.3}(i)
by Lemma \ref{t5.4}(ia) and Proposition \ref{t2.2}.
Similarly, the first assertion to Theorem \ref{tt1.3a}(ii) and
Theorem \ref{tt1.3}(ii)
can be gotten by Lemma \ref{t5.4}(ib) and
Proposition \ref{t3.1b} with the test function $h$ given in \eqref{4.5}.
\qed

\noindent{\it Proof of Theorem \ref{t1.4}.}
Under the assumptions, there are constant $\varepsilon>0$ and $\rho_1,\rho_2>\theta_1\vee\theta_2$
so that \eqref{5.3}--\eqref{5.4} hold.
Let $h$ be the function defined by \eqref{4.5b}
and
 \beqlb\label{4.8b}
g(x,y):=(\delta_0\tilde{c}_1^{-\rho_1}\wedge \tilde{c}_1^{-\rho_2})^\rho-h(x,y),
 \eeqlb
with the constant $\tilde{c}_1>0$ given in Lemma \ref{t5.4}(ii).
Then, using Propositions \ref{t2.2} and \ref{t3.1b} together with Lemma \ref{t5.4}(ii), one obtains the first assertion with test function
$g$, and the second assertion with test function $h$.
\qed

\subsection{Proofs of Theorems \ref{t1.8}--\ref{t1.111}.}

In this subsection we establish the proof of Theorem \ref{t1.8}
by Proposition \ref{t3.1} with the key test function
given in Lemmas \ref{tl3.1}--\ref{t8.5}.
For checking the key condition (ii) in Proposition \ref{t3.1},
we first establish the boundedness of $\mathcal{L}g$
in the following Lemmas \ref{tl3.1}--\ref{t8.5}.
Theorem \ref{t1.9} is proving by using Lemma \ref{t2.3} and
Proposition \ref{t2.4b} with the key test function $g$ chosen in
\eqref{3.6b}.
Theorem \ref{t1.111} can be given by the key corresponding assertions in
Lemma \ref{t5.4}(iii) and Lemmas \ref{t2.3} and \ref{t8.10}.

\blemma\label{tl3.1}
Suppose that $v>0$ and $\theta_i-1<(r_i\vee0)$ for $i=1,2$,
and that one of the conditions (i)--(iii) in Theorem \ref{t1.8} holds.
Let $0\le g\in C^2((0,\infty)\times(0,\infty))$
such that
$g(x,y)=
1+\ln(x^{\rho_1}+y^{\rho_2})+[-\ln(v^{\rho_1}+v^{\rho_2}) ]\vee1$ for $x, y>v$.
Then there are constants $\rho_1,\rho_2,d>0$ such that
$\mathcal{L}g(x,y)\le d g(x,y)$ for all $x,y\ge v$.
\elemma
\proof
We first assume that condition (i) of Theorem \ref{t1.8} holds in the following.
Then there are constants $\rho_1,\rho_2>1$ so that
 \beqnn
\rho_1/\rho_2\le (1-\theta_1)/\kappa_1,\qquad
\rho_2/\rho_1\le (1-\theta_2)/\kappa_2.
 \eeqnn
Then by Lemma \ref{A1.1}(iii), there are constants $C_1,C_2>0$ such that
for all $x,y\ge v$,
 \beqlb\label{2.2}
x^{\theta_1+\rho_1-1}y^{\kappa_1}
\le
(x^{\rho_1})^{1-(1-\theta_1)/\rho_1}(y^\rho_2)^{\kappa_1/\rho_2}
\le C_1(x^{\rho_1}+y^{\rho_2}),~~
y^{\theta_2+\rho_2-1}x^{\kappa_2}\le C_2(x^{\rho_1}+y^{\rho_2}).
 \eeqlb
Observe that
 \beqnn
g_x'(x,y)=\rho_1(x^{\rho_1}+y^{\rho_2})^{-1} x^{\rho_1-1},~~~
g_{xx}''(x,y)
\le
\rho_1(\rho_1-1)(x^{\rho_1}+y^{\rho_2})^{-1}x^{\rho_1-2}
 \eeqnn
and
 \beqnn
g_y'(x,y)=\rho_2 (x^{\rho_1}+y^{\rho_2})^{-1}y^{\rho_2-1},~~~
g_{yy}''(x,y)
\le
\rho_2(\rho_2-1) (x^{\rho_1}+y^{\rho_2})^{-1}y^{\rho_2-2}.
 \eeqnn
Similar to \eqref{4.1},
 \beqlb\label{4.3}
\int_0^\infty K_z^1g(x,y)\mu_1(\dd z)
=x^{-\alpha_1}
\int_0^\infty K(x^{\rho_1}(x^{\rho_1}+y^{\rho_2})^{-1},z;\rho_1)\mu_1(\dd z),
 \eeqlb
where
$
K(v,z;\rho_1):=\ln[v((1+z)^{\rho_1}-1)+1]-\rho_1 zv$ for $v,z,\rho_1>0$.
By Taylor's formula, for $0\le v\le 1$,
 \beqlb\label{4.3aa}
K''_{zz}(v,z;\rho_1)
\le v\rho_1(\rho_1-1)(1+z)^{\rho_1-2}[v((1+z)^{\rho_1}-1)+1]^{-1}
\le \rho_1(\rho_1-1)(1+z)^{-2}.
 \eeqlb
and then by Taylor's formula,
 \beqnn
\int_0^\infty K(v,z;\rho_1)\mu_1(\dd z)
 \ar=\ar
\int_0^\infty z^2\mu_1(\dd z)\int_0^1K''_{zz}(v,zu;\rho_1)(1-u)\dd u \cr
 \ar\le\ar
\rho_1(\rho_1-1)
\int_0^\infty z^2\mu_1(\dd z)\int_0^1(1+zu)^{-2}(1-u)\dd u
=:c_2.
 \eeqnn
It thus follows from \eqref{4.3} and the assumption of $r_{12}-\alpha_1
\le r_1\le0$ that
 \beqnn
x^{r_{12}}\int_0^\infty K_z^1g(x,y)\mu_1(\dd z)
\le c_2 x^{r_{12}-\alpha_1}\le
c_2 v^{r_{12}-\alpha_1},\quad x\ge v.
 \eeqnn
Similarly, by the assumption of $r_{22}-\alpha_2
\le r_2\le0$,
 \beqlb\label{4.3aaa}
y^{r_{22}}\int_0^\infty K_z^2g(x,y)\mu_2(\dd z)
\le c_3  v^{r_{22}-\alpha_2},\quad y\ge v
 \eeqlb
for some constant $c_3>0$.
Then, by \eqref{3.0} and \eqref{2.2},
there are constants $C_3,C_4,C_5>0$ such that
 \beqnn
 \ar\ar
(x^{\rho_1}+y^{\rho_2}) \mathcal{L}g(x,y) \cr
 \ar\le\ar
a_1\rho_1x^{\theta_1+\rho_1-1}y^{\kappa_1}
+a_2\rho_2y^{\theta_2+\rho_2-1}x^{\kappa_2}
+C_3(x^{\rho_1+r_{11}-2}+y^{\rho_2+r_{21}-2})
+c_2 v^{r_{12}-\alpha_1}
+c_3 v^{r_{22}-\alpha_2} \cr
 \ar\le\ar
C_4(x^{\rho_1}+y^{\rho_2}),\qquad x,y\ge v,
 \eeqnn
where the last inequality follows from
the assumptions $r_{11}-2,r_{21}-2\le 0$.
This gives the assertion.

Since the proofs are similar under conditions (ii) and (iii)  of Theorem \ref{t1.8}, we only state that under condition (iii).
Then $\frac{\kappa_1}{r_1+1-\theta_1}\le\frac{1-\theta_2}{\kappa_2}$.
Thus there are constants $1-\theta_1<\rho_1<1$ and $\rho_2>1-\theta_2$ so that
$\frac{\kappa_1}{r_1+1-\theta_1}\le\frac{\rho_2}{r_1+\rho_1}
\le\frac{1-\theta_2}{\kappa_2}$,
which implies
 \beqlb\label{4.4}
\frac{\theta_1+\rho_1-1}{r_1+\rho_1}+\frac{\kappa_1}{\rho_2}\le1,
\qquad
\frac{\theta_2+\rho_2-1}{\rho_2}+\frac{\kappa_2}{r_1+\rho_1}\le1.
 \eeqlb
By the first inequality of \eqref{4.3aa}, for $0<\rho_1<1$ and $0\le v\le 1$,
 \beqnn
K''_{zz}(v,z;\rho_1)
\le v\rho_1(\rho_1-1)(1+z)^{\rho_1-2}[v((1+z)^{\rho_1}-1)+1]^{-1}
\le v\rho_1(\rho_1-1)(1+z)^{-2}.
 \eeqnn
and then by Taylor's formula,
 \beqnn
\int_0^\infty K(v,z;\rho_1)\mu_1(\dd z)
 \ar=\ar
\int_0^\infty z^2\mu_1(\dd z)\int_0^1K''_{zz}(v,zu;\rho_1)(1-u)\dd u \cr
 \ar\le\ar
\rho_1(\rho_1-1)v
\int_0^\infty z^2\mu_1(\dd z)\int_0^1(1+zu)^{-2}(1-u)\dd u
=:\rho_1(\rho_1-1)c_4v.
 \eeqnn
It thus follows from \eqref{4.3} that
 \beqnn
\int_0^\infty K_z^1g(x,y)\mu_1(\dd z)
\le c_4\rho_1(\rho_1-1)(x^{\rho_1}+y^{\rho_2})^{-1}x^{\rho_1-\alpha_1}.
 \eeqnn
Now by \eqref{3.0} and \eqref{4.3aaa}, $x,y\ge v$,
 \beqlb\label{4.6}
 \ar\ar
[x^{\rho_1}+y^{\rho_2}] \mathcal{L}g(x,y) \cr
 \ar\le\ar
a_1\rho_1x^{\theta_1+\rho_1-1}y^{\kappa_1}
+a_2\rho_2y^{\theta_2+\rho_2-1}x^{\kappa_2}
-\rho_1 b_{10}x^{\rho_1+r_{10}-1}  \cr
 \ar\ar
-b_{11} \rho_1(1-\rho_1)x^{\rho_1+r_{11}-2}
-b_{12} \rho_1(1-\rho_1)c_4x^{\rho_1+r_{12}-\alpha_1}
+b_{20}\rho_2y^{\rho_2+r_{20}-1} \cr
 \ar\ar
+b_{21} \rho_2(\rho_2-1)1_{\{\rho_2>1\}}y^{\rho_2+r_{21}-2}
+b_{22} c_31_{\{\rho_2>1\}}v^{r_{22}-\alpha_2}[x^{\rho_1}+y^{\rho_2}] \cr
 \ar\le\ar
a_1\rho_1x^{\theta_1+\rho_1-1}y^{\kappa_1}
+a_2\rho_2y^{\theta_2+\rho_2-1}x^{\kappa_2}
-\tilde{b}_1x^{\rho_1+r_1}
+\tilde{b}_2[x^{\rho_1}+y^{\rho_2}]
 \eeqlb
for some constants $\tilde{b}_1,\tilde{b}_2>0$.
By \eqref{4.4} and Lemma \ref{A1.1}(iii) again, there is a constant $C_1>0$ such that
 \beqnn
a_1\rho_1x^{\theta_1+\rho_1-1}y^{\kappa_1}
+a_2\rho_2y^{\theta_2+\rho_2-1}x^{\kappa_2}\le \tilde{b}_1x^{r_1+\rho_1}+C_1y^{\rho_2},
\qquad x,y\ge v.
 \eeqnn
Combining this with \eqref{4.6} one concludes the assertion.
\qed

\blemma\label{t8.2}
Let $v>0$ be fixed.
Suppose that $(\theta_i-1)\vee0<r_i$ holds for $i=1,2$
and the assumptions in Theorem \ref{t1.8}(iva) hold,
and that
$g(x,y):=(x^{\rho_1}+y^{\rho_2})^\rho$ for $x,y>0$.
Then,
there are constants $\rho_1,\rho_2,d>0$ and $0<\rho<1$ such that
$\mathcal{L}g(x,y)\le d g(x,y)$ for all $x,y\ge v$.
\elemma
\proof
Since the proofs are similar, we only show that
 under Condition \ref{c3}(i).
Since \eqref{1.7a} and \eqref{1.7} hold,
then there are large enough constants $\rho_1,\rho_2\ge2
\vee\kappa_2\vee\kappa_1$ satisfying \eqref{8.1}
and
 \beqlb\label{8.1a}
r_{12}-\alpha_1\le \rho_1 r_2/\rho_2 \mbox{ if } b_{12}>0,\qquad
r_{22}-\alpha_2\le \rho_2 r_1/\rho_1 \mbox{ if } b_{22}>0.
 \eeqlb
Put $0<\rho<1$ and
 \beqnn
J_1(x,y):=
b_{10}\rho_1x^{r_{10}-1+\rho_1}
+b_{20}\rho_2y^{r_{20}-1+\rho_2},~~
J_2(x,y):= a_1\rho_1x^{\theta_1-1+\rho_1}y^{\kappa_1}
+a_2\rho_2 x^{\kappa_2}y^{\theta_2-1+\rho_2}.
 \eeqnn
Under Condition \ref{c3}(i),
$J_1(x,y)=
b_1\rho_1x^{r_1+\rho_1}
+b_2\rho_2y^{r_2+\rho_2}$.
Let $K(v,z)$ be the function defined by \eqref{3.7bb}.
Then by the arguments in \eqref{4.1} and Lemma \ref{t8.1}(iiia) we obtain
 \beqnn
 \ar\ar
\int_0^\infty K_z^1g(x,y)\mu_1(\dd z)
=x^{-\alpha_1}\int_0^\infty K_{xz}^1g(x,y)\mu_1(\dd z) \cr
 \ar=\ar
-x^{-\alpha_1}(x^{\rho_1}+y^{\rho_2})^\rho
\int_0^\infty
K(x^{\rho_1}(x^{\rho_1}+y^{\rho_2})^{-1},z)\mu_1(\dd z)
\le
\rho \tilde{d} (x^{\rho_1}+y^{\rho_2})^{\rho-\alpha_1/\rho_1}
 \eeqnn
and
 \beqnn
\int_0^\infty K_z^2g(x,y)\mu_2(\dd z)
\le
\rho \tilde{d} (x^{\rho_1}+y^{\rho_2})^{\rho-\alpha_2/\rho_2}
 \eeqnn
for some constant $\tilde{d}>0$.
Let
 \beqnn
J_3(x,y)
:=
\rho_1(\rho_1-1) b_{11}x^{r_{11}+\rho_1-2}
+\rho_2(\rho_2-1) b_{21}y^{r_{21}+\rho_2-2}.
 \eeqnn
Then by \eqref{3.0},
 \beqlb\label{3.4aaa}
\mathcal{L}g(x,y)
\le
\rho[(x^{\rho_1}+y^{\rho_2})]^{\rho-1}[-J_1(x,y)+J_2(x,y)+J_3(x,y)
+\tilde{d}J_4(x,y)],
 \eeqlb
where
 \beqnn
J_4(x,y):=b_{12}x^{r_{12}}(x^{\rho_1}+y^{\rho_2})^{1-\alpha_1/\rho_1}
+b_{22}y^{r_{22}}(x^{\rho_1}+y^{\rho_2})^{1-\alpha_2/\rho_2}.
 \eeqnn 
By \eqref{8.1} and Lemma \ref{t8.1}(i),
there is a constant $c_1>1$ so that
$4^{-1}J_1(x,y)\ge J_2(x,y)$ for all
$x,y\ge c_1$.
There is a constant $c_2>c_1$ so that
 \beqnn
8^{-1}J_1(x,y)-J_2(x,y)
\ge
8^{-1}b_1\rho_1x^{r_1+\rho_1}
-a_1\rho_1x^{\theta_1-1+\rho_1} c_1^{\kappa_1}
-a_2\rho_2 c_1^{\theta_2-1+\rho_2}x^{\kappa_2}\ge0
 \eeqnn
for all $x\ge c_2$ and $0<y\le c_1$, and
 \beqnn
8^{-1}J_1(x,y)-J_2(x,y)
\ge
8^{-1}b_2\rho_2y^{r_2+\rho_2}
-a_1\rho_1 c_1^{\theta_1-1+\rho_1}y^{\kappa_1}
-a_2\rho_2 y^{\theta_2-1+\rho_2}c_1^{\kappa_2}\ge0
 \eeqnn
for all $y\ge c_2$ and $0<x\le c_1$.
Thus,
 \beqlb\label{3.5}
8^{-1}J_1(x,y)\ge J_2(x,y),\qquad x\ge c_2\mbox{ or }y\ge c_2.
 \eeqlb
Under Condition \ref{c3}(i), there is a constant $c_3>c_2$ such that
 \beqlb\label{3.5b}
8^{-1}J_1(x,y)\ge J_3(x,y),\qquad x\ge c_3\mbox{ or }y\ge c_3.
 \eeqlb
By Lemma \ref{A1.1}(ii),
 \beqnn
J_4(x,y)\le
b_{12}x^{r_{12}+\rho_1-\alpha_1}
+b_{22}y^{r_{22}+\rho_2-\alpha_2}
+b_{12}x^{r_{12}} y^{\rho_2-\rho_2\alpha_1/\rho_1}
+b_{22}y^{r_{22}}x^{\rho_1-\rho_1\alpha_2/\rho_2}.
 \eeqnn
By \eqref{8.1a},
 \beqlb\label{3.5aa}
\frac{r_{12}}{r_1+\rho_1}\le \frac{r_2+\rho_2\alpha_1/\rho_1}{r_2+\rho_2}
\mbox{ if } b_{12}>0,\qquad
\frac{r_{22}}{r_2+\rho_2}\le \frac{r_1+\rho_1\alpha_2/\rho_2}{r_1+\rho_1} \mbox{ if } b_{22}>0
 \eeqlb
and then
 \beqlb\label{3.5a}
\frac{r_{12}}{r_1+\rho_1}+\frac{\rho_2-\rho_2\alpha_1/\rho_1}{r_2+\rho_2}\le1
 \mbox{ if } b_{12}>0,\qquad
\frac{r_{22}}{r_2+\rho_2}+\frac{\rho_1-\rho_1\alpha_2/\rho_2}{r_1+\rho_1}\le1
 \mbox{ if } b_{22}>0
 \eeqlb
Then by \eqref{3.5a}, similar to \eqref{3.5} and \eqref{3.5b},
there is a constant $c_4>c_3$ such that
 \beqlb\label{3.5c}
8^{-1}J_1(x,y)\ge \tilde{d}J_4(x,y),\qquad x\ge c_4\mbox{ or }y\ge c_4.
 \eeqlb
Combining \eqref{3.5}--\eqref{3.5b} and \eqref{3.5c} with \eqref{3.4aaa} we obtain
$\mathcal{L}g(x,y)\le 0$ whenever $x\ge c_4$ or $y\ge c_4$.
Observe that there is a constant $C_1>0$ such that for all
$v\le x,y\le c_4$, we have
$J_2(x,y)+J_3(x,y)+\tilde{d}J_4(x,y)\le C_1$, and then by \eqref{3.4aaa}
 \beqnn
\mathcal{L}g(x,y)
\le
\rho[(x^{\rho_1}+y^{\rho_2})]^{\rho-1}[J_2(x,y)+J_3(x,y)+\tilde{d}J_4(x,y)]
\le
C_1\rho [(v^{\rho_1}+v^{\rho_2})]^{\rho-1},v\le x,y\le c_4.
 \eeqnn
Since $g$ are bounded on $[v,c_4]\times[v,c_4]$,
the proof is concluded.
\qed

\blemma\label{t8.3}
Let $v>0$.
Suppose that $(\theta_i-1)\vee0<r_i$ holds for $i=1,2$
and the assumptions in Theorem \ref{t1.8}(ivb) hold,
and that
$g(x,y):=(\delta_0x^{\rho_1}+y^{\rho_2})^\rho$.
Then,
there are constants $\rho_1,\rho_2,\delta_0,d>0$ and $0<\rho<1$
such that
$\mathcal{L}g(x,y)\le d g(x,y)$ for all $x,y\ge v$.
\elemma
\proof
The proof is a modification of that of Lemma \ref{t8.2},
we omit it here.
\qed

\blemma\label{t8.5}
Suppose that $(\theta_i-1)\vee0<r_i$ holds for $i=1,2$, and the assumptions in Theorem \ref{t1.8}(ivc) hold.
Let
$g(x,y):=\tilde{h}(x,y)^\rho$
with
 \beqlb\label{8.23a}
\tilde{h}(x,y):=x^{\rho_1}+ y(1+y)^{-\varepsilon},\qquad x,y>0
 \eeqlb
when \eqref{1.5} holds,
and
 \beqlb\label{8.23b}
\tilde{h}(x,y):=y^{\rho_2}+ x(1+x)^{-\varepsilon},\qquad x,y>0
 \eeqlb
when \eqref{1.6} holds.
Then
there are constants $0<\rho<1$, $0<\varepsilon<1/2$ and $\rho_1,\rho_2,d>0$ such that
$\mathcal{L}g(x,y)\le d g(x,y)$ for all $x,y\ge v$.
\elemma
\proof
Since the proofs are similar, we only present that under \eqref{1.5}.
Under \eqref{1.7} and \eqref{1.5},
there is a constant $\rho_1>1$ so that
 \beqnn
1<\rho_1,\quad r_1<\rho_1r_2,\quad
\frac{\kappa_2(r_2+1)}{r_2+1-\theta_2}-r_1
<\rho_1
<\frac{(r_1+1-\theta_1)(r_2+1)}{\kappa_1}-r_1.
 \eeqnn
Then
 \beqnn
\frac{r_1+\rho_1}{r_1+2\rho_1}
+\frac{1}{r_2+2}<1,~
\frac{\rho_1}{r_1+\rho_1}>
\frac{1}{r_2+1},~
\frac{\theta_1-1+\rho_1}{r_1+\rho_1}
+\frac{\kappa_1}{r_2+1}<1,~
\frac{\kappa_2}{r_1+\rho_1}
+\frac{\theta_2}{r_2+1}<1.
 \eeqnn
Thus, there is a small enough constant $0<\varepsilon<\frac{1}{2}\wedge r_2$ so that
$r_2+1-\theta_2>2\varepsilon$, $\frac{\rho_1}{r_1+\rho_1}>
\frac{1-\varepsilon}{r_2+1-\varepsilon}$,
 \beqlb\label{8.9}
\frac{r_1+\rho_1}{r_1+2\rho_1}
+\frac{1}{r_2+2-2\varepsilon}<1
 \eeqlb
and
 \beqlb\label{8.10}
\frac{\theta_1-1+\rho_1}{r_1+\rho_1}
+\frac{\kappa_1}{r_2+1-\varepsilon}<1,\qquad
\frac{\kappa_2}{r_1+\rho_1}
+\frac{\theta_2}{r_2+1-\varepsilon}<1.
 \eeqlb
For $0<\rho<1/\rho_1$, we have
 \beqlb\label{8.9b}
\frac{r_1+\rho\rho_1}{r_1+\rho_1}
+\frac{(1-\varepsilon)(1-\rho)}{r_2+1-\varepsilon}
<1.
 \eeqlb

It is elementary to see that
$g_x'(x,y)=\rho\rho_1x^{\rho_1-1}\tilde{h}(x,y)^{\rho-1}$
and
 \beqnn
 \ar\ar
g_{xx}''(x,y) \cr
 \ar=\ar
-\rho\rho_1^2(1-\rho)x^{2\rho_1-2}\tilde{h}(x,y)^{\rho-2}
+\rho\rho_1(\rho_1-1+\delta) x^{\rho_1-2}\tilde{h}(x,y)^{\rho-1}
-\rho\rho_1\delta x^{\rho_1-2}\tilde{h}(x,y)^{\rho-1} \cr
 \ar\le\ar
-\rho\rho_1\delta x^{\rho_1-2}\tilde{h}(x,y)^{\rho-1}
-\rho\rho_1\tilde{h}(x,y)^{\rho-2}[(1-\rho_1\rho-\delta)x^{2\rho_1-2}
-(\rho_1-1+\delta)x^{\rho_1-2}y]
 \eeqnn
for $0<\delta<(1-\rho_1\rho)[c_1(\rho\rho_1)\wedge1]$ with
	$c_1(\rho\rho_1)$ defined by \eqref{3.0a}.

By the arguments in \eqref{4.1},
 \beqnn
\int_0^\infty K_z^1g(x,y)\mu_1(\dd z)
=
-\tilde{h}(x,y)^\rho x^{-\alpha_1}\int_0^\infty K(\frac{x^{\rho_1}}{\tilde{h}(x,y)},z)\mu_1(\dd z)
 \eeqnn
In view of Lemma \ref{t8.1}(iiib),
 \beqnn
 \ar\ar
\int_0^\infty K_z^1g(x,y)\mu_1(\dd z) \cr
 \ar\le\ar
-\rho\rho_1  x^{-\alpha_1}\Big[(1-\rho\rho_1)c_1(\rho\rho_1)
x^{2\rho_1}\tilde{h}(x,y)^{\rho-2}
-(\rho_1-1)d_{1,\sigma}
x^{\rho_1}y\tilde{h}(x,y)^{\rho-2}
-(\rho_1-1)d_{2,\sigma}
x^{\rho\rho_1} \Big] \cr
 \ar\le\ar
d_{2,\sigma}\rho\rho_1(\rho_1-1) x^{\rho\rho_1-\alpha_1}
-\rho\rho_1\delta x^{\rho_1-\alpha_1}\tilde{h}(x,y)^{\rho-1} \cr
 \ar\ar
-\rho\rho_1\tilde{h}(x,y)^{\rho-2}\Big[[(1-\rho_1\rho)
c_1(\rho_1\rho)-\delta]x^{2\rho_1-\alpha_1}
-[(\rho_1-1)d_{1,\sigma}+\delta]x^{\rho_1-\alpha_1}y\Big],
 \eeqnn
where $K(v,z)$ is defined in \eqref{3.7bb} and
$\lim_{\sigma\to\infty}d_{2,\sigma}=0$.
Similarly,
 \beqnn
g_y'(x,y)
=
\rho \tilde{h}(x,y)^{\rho-1}
(1+y)^{-1-\varepsilon}[1+(1-\varepsilon )y]
 \eeqnn
and
 \beqnn
g_y'(x,y)
\ge \rho (1-\varepsilon)\tilde{h}(x,y)^{\rho-1}
(1+y)^{-\varepsilon}
\ge 2^{-1}\rho \tilde{h}(x,y)^{\rho-1}
(1+y)^{-\varepsilon}.
 \eeqnn
Moreover,
 \beqnn
g_{yy}''(x,y)
 \le
-2^{-2}\rho(1-\rho)\tilde{h}(x,y)^{\rho-2}
(1+y)^{-2\varepsilon}
-2^{-1}\varepsilon\rho \tilde{h}(x,y)^{\rho-1}
(1+y)^{-1-\varepsilon}.
 \eeqnn
For all $y>0$ and $0<z\le1$, we have
$\tilde{h}(x,y(1+z))\le 2\tilde{h}(x,y)$
and then
 \beqnn
g_{yy}''(x,y(1+z))
 \le
-2^{\rho-2\varepsilon-4}\rho(1-\rho)\tilde{h}(x,y)^{\rho-2}
(1+y)^{-2\varepsilon}
-2^{\rho-3-\varepsilon}\varepsilon\rho \tilde{h}(x,y)^{\rho-1}
(1+y)^{-1-\varepsilon}.
 \eeqnn
From Taylor's formula and replacing the variable $z$ by $yz$ it follows that
 \beqnn
 \ar\ar
\int_0^\infty K_z^2g(x,y)\mu_2(\dd z)
=\int_0^\infty z^2\mu_2(\dd z)\int_0^1g_{yy}''(x,y+zu)(1-u)\dd u \cr
 \ar\le\ar
y^{2-\alpha_2}
\int_0^1 z^2\mu_2(\dd z)\int_0^1g_{yy}''(x,y+yzu)(1-u)\dd u\cr
 \ar\le\ar
-\rho y^{2-\alpha_2}
\Big[\tilde{c}_1(1-\rho)\tilde{h}(x,y)^{\rho-2}(1+y)^{-2\varepsilon}
+\varepsilon \tilde{c}_2\tilde{h}(x,y)^{\rho-1}(1+y)^{-1-\varepsilon}\Big]
 \eeqnn
for some constants $\tilde{c}_1,\tilde{c}_2>0$.

Now by \eqref{3.0} and \eqref{1.6b},
 \beqlb\label{8.8}
\mathcal{L}g(x,y)
\le
-\rho \tilde{h}(x,y)^{\rho-1} [I_2(x,y)-I_1(x,y)]
-\rho \tilde{h}(x,y)^{\rho-2} I_3(x,y)
+d_{2,\sigma}\rho I_4(x,y),
 \eeqlb
where $I_1(x,y):=\rho_1a_1x^{\theta_1-1+\rho_1}y^{\kappa_1}
+a_2y^{\theta_2}x^{\kappa_2}$,
$I_4(x,y):=\rho_1(\rho_1-1)b_{12} x^{r_{12}-\alpha_1+\rho\rho_1}$,
 \beqnn
I_2(x,y):=\rho_1\delta   b_1x^{r_1+\rho_1}
+2^{-1}(\tilde{c}_2\wedge1)\varepsilon  b_2y^{2+r_2}(1+y)^{-1-\varepsilon}
 \eeqnn
and
 \beqnn
I_3(x,y)
 \ar:=\ar
\rho_1[(1-\rho_1\rho)(c_1(\rho_1\rho)\wedge1)-\delta] b_1x^{r_1+2\rho_1}
+ 2^{-2}(1-\rho)(\tilde{c}_1\wedge1) b_2y^{2+r_2}(1+y)^{-2\varepsilon} \cr
 \ar\ar
-\rho_1\big[(\rho_1-1+\delta)b_{11}+
[(\rho_1-1)d_{1,\sigma}+\delta]b_{12}\big]x^{r_1+\rho_1}y.
 \eeqnn

By Lemma \ref{A1.1}(i) and \eqref{8.10},
there is a constant $c_0>1$ so that
for all $x,y\ge c_0$,
 \beqlb\label{8.12}
I_2(x,y)
\ge
\rho_1\delta   b_1x^{r_1+\rho_1}
+2^{-2-\varepsilon}(\tilde{c}_2\wedge1)\varepsilon b_2y^{2+r_2}y^{-1-\varepsilon}
\ge 4I_1(x,y),\qquad x,y\ge c_0.
 \eeqlb
By Lemma \ref{A1.1}(i) and \eqref{8.9}, there is a constant $c_1>c_0$ so that
$I_3(x,y)\ge0$ for all $x,y\ge c_1$.
By Lemma \ref{A1.1}(ii),
 \beqnn
x^{r_{12}-\alpha_1+\rho\rho_1}\tilde{h}(x,y)^{1-\rho}
\le
x^{r_{12}-\alpha_1+\rho_1}
+x^{r_{12}-\alpha_1+\rho\rho_1}
y^{(1-\varepsilon)(1-\rho)},\qquad x,y>0.
 \eeqnn
By Lemma \ref{A1.1}(i) and \eqref{8.9b},
there is a constant $c_2>c_1$ so that
 \beqnn
I_2(x,y)
\ge 8d_{2,\sigma}\rho_1(\rho_1-1)b_{12} x^{r_{12}-\alpha_1+\rho\rho_1}
y^{(1-\varepsilon)(1-\rho)},
\qquad x,y\ge c_2
 \eeqnn
and
then for large enough $\sigma>0$,
 \beqnn
I_2(x,y)
\ge
4d_{2,\sigma}I_4(x,y)
\tilde{h}(x,y)^{1-\rho},\qquad x,y\ge c_2,
 \eeqnn
where we use the fact $\lim_{\sigma\to\infty}d_{2,\sigma}=0$.
Since $\theta_1-1<r_1$ and
$\varepsilon+\theta_2-1<r_2$ under the assumptions
and $\kappa_1<r_2+1-\varepsilon$
and $\kappa_2<r_1+\rho_1$ by \eqref{8.10},
then there is a constant $c_3>c_2$ so that
 \beqlb\label{8.3}
I_2(x,y)
\ge
\rho_1\delta   b_1x^{r_1+\rho_1}
\ge
4a_1\rho_1x^{\theta_1-1+\rho_1}  c_2^{\kappa_1}
+4a_2 x^{\kappa_2} c_2^{\theta_2}\ge 4I_1(x,y)
 \eeqlb
for $x\ge c_3$ and $0<y\le c_2$ and
 \beqlb\label{8.4}
I_2(x,y)
\ge
2^{-2-\varepsilon}(\tilde{c}_2\wedge1)\varepsilon b_2y^{2+r_2}y^{-1-\varepsilon}
\ge
4a_1\rho_1y^{\kappa_1} c_2^{\theta_1-1+\rho_1}
+4a_2 y^{\theta_2}  c_2^{\kappa_2}\ge 4I_1(x,y)
 \eeqlb
for $y\ge c_3$ and $0<x\le c_2$.
Similarly, there is a constant $c_4>c_3$
so that $I_3(x,y)\ge 0$ and
$I_2(x,y)
\ge
4d_{2,\sigma}I_4(x,y)
\tilde{h}(x,y)^{1-\rho}$ when $x\ge c_4$ and $0<y\le c_2$, or
$y\ge c_4$ and $0<x\le  c_2$.
Combining these with \eqref{8.8} one gets
$\mathcal{L}g(x,y)\le0$ when $x\ge c_4$ or $y\ge c_4$.
The assertion then follows by an argument identical to that given at the conclusion of the proof of Lemma \ref{t8.2}.
\qed

Now we can get the proof of Theorem \ref{t1.8}.

\noindent{\it Proof of Theorem \ref{t1.8}.}
Let $g$ be the function defined in
Lemmas \ref{tl3.1}--\ref{t8.5} for different conditions.
One concludes the assertion of Theorem \ref{t1.8} by Proposition \ref{t3.1} and Lemmas \ref{tl3.1}--\ref{t8.5}.
\qed

Before proving Theorem \ref{t1.9}, we first state the following lemma.

\blemma\label{t2.3}
Suppose that $(\theta_i-1)\vee0< r_i$ for $i=1,2$.
Let $h(x,y):=x^{\rho_1}+y^{\rho_2}$ under condition (iva) in Theorem \ref{t1.8},
$h(x,y):=\delta_0x^{\rho_1}+y^{\rho_2}$ under condition (ivb) in Theorem \ref{t1.8},
$h(x,y)$ be defined in Lemma \ref{t8.5} under condition (ivc) of Theorem \ref{t1.8}.
We also assume that $h(x,y):=\delta_0x^{\rho_1}+y^{\rho_2}$ when Condition \ref{c3}(i), \eqref{6.2} and the assumptions in Theorem \ref{t1.111}(ii) are satisfied.
For $\rho>0$ let
 \beqlb\label{3.6b}
g(x,y)=1+ h(x,y)^{-\rho},\qquad x,y>0.
 \eeqlb
Then, there are constants $\rho_1,\rho_2,\varepsilon,\rho,\tilde{\rho},\delta_0,C>0$
and $c>1$ so that
$\mathcal{L}g(x,y)\ge 2C a^{\tilde{\rho}}
\ge C a^{\tilde{\rho}}g(x,y)$ for all $x,y>0$ satisfying $h(x,y)\ge a$
and $a>c$.
\elemma
\proof
Since the proofs are similar, then we give the details
to that under condition (iva) in Theorem \ref{t1.8}
and Condition \ref{c3}(i), which is a modification to that of Lemma \ref{t8.2}.
Let $0<\rho <(r_1/\rho_1)\wedge (r_2/\rho_2)$.
Similar to \eqref{3.4aaa},  we obtain
 \beqlb\label{3.5bbda}
\mathcal{L}g(x,y)
\ge
\rho[(x^{\rho_1}+y^{\rho_2})]^{-\rho-1}[J_1(x,y)-J_2(x,y)-J_3(x,y)
-\tilde{d}J_4(x,y)].
 \eeqlb
where $\tilde{d}>0$ and $J_i(x,y)$ for $i=1,2,3,4$ are defined in the proof of Lemma \ref{t8.2}.
By \eqref{3.5} and  \eqref{3.5b},
there is a constant $c_1>0$ such that
 \beqnn
2^{-1}J_1(x,y)\ge J_2(x,y)+J_3(x,y)+\tilde{d}J_4(x,y) \mbox{ for all }x,y>0 \mbox{ satisfying }
x^{\rho_1}+y^{\rho_2}\ge c_1.
 \eeqnn
For $a>c_1/2$ and $x^{\rho_1}+y^{\rho_2}\ge 2a$,
when $x^{\rho_1}\le y^{\rho_2}$, we have $y^{\rho_2}\ge a$ and then
 \beqnn
[(x^{\rho_1}+y^{\rho_2})]^{-\rho-1}J_1(x,y)
\ge 2^{-\rho-1}b_2\rho_2y^{-\rho_2(\rho+1)}y^{r_2+\rho_2}
=2^{-\rho-1}b_2\rho_2 y^{r_2-\rho_2\rho}
\ge 2^{-\rho-1}b_2 a^{(r_2-\rho_2\rho)/\rho_2}.
 \eeqnn
Similarly, for $a>c_1/2$ and $x^{\rho_1}+y^{\rho_2}\ge 2a$,
when $x^{\rho_1}> y^{\rho_2}$,
we have $x^{\rho_1}\ge a$ and then
$[(x^{\rho_1}+y^{\rho_2})]^{-\rho-1}J_1(x,y)
\ge 2^{-\rho-1}b_1 a^{(r_1-\rho_1\rho)/\rho_1}$.
Thus for all $a>c_1$,
 \beqlb\label{3.5bbd}
[(x^{\rho_1}+y^{\rho_2})]^{-\rho-1}J_1(x,y)
\ge 2^{-\rho-1}(b_1 \wedge  b_2)
a^{\tilde{\rho}},~
\mbox{ for all } x,y>0 \mbox{ satisfying }x^{\rho_1}+y^{\rho_2}\ge 2a
 \eeqlb
with $\tilde{\rho}:=(r_1/\rho_1)\wedge (r_2/\rho_2)-\rho>0$.
Combining \eqref{3.5bbda} with  \eqref{3.5bbd} we obtain
 \beqlb\label{4.10}
\mathcal{L}g(x,y)
\ge
2^{-\rho-2}\rho(b_1\wedge b_2)a^{\tilde{\rho}},\qquad
\mbox{ for all } x,y>0 \mbox{ satisfying }x^{\rho_1}+y^{\rho_2}\ge 2a,
 \eeqlb
which gives the assertion.

\noindent{\it Proof of Theorem \ref{t1.9}.}
By Lemma \ref{t2.3} and Proposition \ref{t2.4b}, one gets the
assertion immediately.
\qed

Before stating the proof of Theorem \ref{t1.111}, we show the following
key assertion.

\blemma\label{t8.10}
Suppose that $(\theta_i-1)\vee0< r_i$ for $i=1,2$, and Condition \ref{c3}(i)
and \eqref{6.2} are satisfied.
Then there are constants $\rho\in(0,1)$ and $\rho_1,\rho_2,\delta_0,d>0$ such that
$\mathcal{L}g(x,y)\le dg(x,y)$ for all $x,y\ge v$ if
the assumptions in Theorem \ref{t1.111}(iii) holds,
where $g(x,y):=(\delta_0x^{\rho_1}+y^{\rho_2})^{\rho}$ for $x,y>0$.
\elemma
\proof
We first consider the assertion under the case of $\bar{r}_i<0$ for $i=1,2$.
Under \eqref{6.2},
there are constants $\rho_1,\rho_2>1$ such that
 \beqlb\label{1.12}
\frac{r_1+1-\theta_1}{\kappa_1}=\frac{r_1+\rho_1}{r_2+\rho_2}
=\frac{\kappa_2}{ r_2+1-\theta_2}.
 \eeqlb
Then, using \eqref{1.12} and \cite[Lemma 3.8]{XYZh25}, there is a constant $\delta_0>0$
such that
 \beqlb\label{6.1}
J_2(x,y)\le J_1(x,y),\qquad x,y>0,
 \eeqlb
where
 \beqnn
J_1(x,y):=
b_{10}\rho_1\delta_0x^{r_{10}-1+\rho_1}
+b_{20}\rho_2y^{r_{20}-1+\rho_2},
J_2(x,y):= a_1\rho_1\delta_0x^{\theta_1-1+\rho_1}y^{\kappa_1}
+a_2\rho_2 x^{\kappa_2}y^{\theta_2-1+\rho_2}.
 \eeqnn
Let $\rho\in(0,1)$ be small enough such that
 \beqlb\label{3.4aaac}
\bar{r}_1+\rho\rho_1<0,\qquad
\bar{r}_2+\rho\rho_2<0.
 \eeqlb
Similar to \eqref{3.4aaa},
we have
 \beqlb\label{3.4aaab}
\mathcal{L}g(x,y)
\le
\rho[\delta_0x^{\rho_1}+y^{\rho_2}]^{\rho-1}
[-J_1(x,y)+J_2(x,y)+J_3(x,y)+\tilde{d}J_4(x,y)]
 \eeqlb
for some constant $\tilde{d}>0$,
where
 \beqnn
J_3(x,y)
 \ar:=\ar
\rho_1(\rho_1-1) \delta_0b_{11}x^{r_{11}+\rho_1-2}
+\rho_2(\rho_2-1) b_{21}y^{r_{21}+\rho_2-2}
 \eeqnn
and
 \beqnn
J_4(x,y):=\delta_0b_{12}x^{r_{12}}(x^{\rho_1}+y^{\rho_2})^{1-\alpha_1/\rho_1}
+b_{22}y^{r_{22}}(x^{\rho_1}+y^{\rho_2})^{1-\alpha_2/\rho_2}.
 \eeqnn
Since $\bar{r}_i<0$ for $i=1,2$, we have \eqref{8.1a} and \eqref{3.5a}
with $r_i$ replaced by $r_{i2}$.
Then similar to \eqref{3.5c}, we have $J_4(x,y)\le C \tilde{J}_4(x,y):= b_{11}x^{r_{12}+\rho_1-\alpha_1}
+b_{21}y^{r_{22}+\rho_2-\alpha_2}$
for all $x,y\ge v$ and some constant $C>0$.
By \eqref{3.4aaac}, there are constants $C_1,C_2,C_3>0$ such that
 \beqnn
J_3(x,y)+C\tilde{d} \tilde{J}_4(x,y)
 \ar\le\ar
C_1x^{\bar{r}_1+\rho_1}
+C_2y^{\bar{r}_2+\rho_2}
=C_1x^{\bar{r}_1+\rho\rho_1}x^{(1-\rho)\rho_1}
+C_2y^{\bar{r}_2+\rho\rho_2}y^{(1-\rho)\rho_2} \cr
 \ar\le\ar
C_1v^{\bar{r}_1+\rho\rho_1}x^{(1-\rho)\rho_1}
+C_2v^{\bar{r}_2+\rho\rho_2}y^{(1-\rho)\rho_2}
\le
C_3(\delta_0x^{\rho_1})^{1-\rho}
+C_3y^{(1-\rho)\rho_2} \cr
 \ar\le\ar
2C_3(\delta_0x^{\rho_1}+y^{\rho_2})^{1-\rho},\qquad x,y\ge v,
 \eeqnn
which implies
$\mathcal{L}g(x,y)\le 2\rho C_3$ for all $x,y\ge v$.
This gives the assertion since $\inf_{x,y\ge v}g(x,y)>0$.

In the following we show the assertion under \eqref{1.11} and we only state that
for the first case of \eqref{1.11}.
Let $M=(r_1+1-\theta_1)/\kappa_1$
and $f(x)=\frac{r_1+x}{r_2+1-x}$ for $x\in[0,1]$.
Then, under the first assumption of \eqref{1.11},
$f(1-\theta_1)\le M\le f(\theta_2)$, and there is a constant $\lambda\in[1-\theta_1,\theta_2]$
such that $f(\lambda)=M$.
Taking $\rho_1=\lambda\in[1-\theta_1,1]$ and $\rho_2=1-\lambda\in[1-\theta_2,1]$ we get
 \beqnn
M=f(\lambda)=\frac{r_1+\lambda}{r_2+1-\lambda}
=\frac{r_1+\rho_1}{r_2+\rho_2},
 \eeqnn
which gives \eqref{1.12}.
Then by the arguments in \cite[Lemma 3.8]{XYZh25},
we obtain \eqref{6.1}.
By Lemma \ref{t8.1}(iiic), similar to \eqref{3.4aaa}, we obtain
$\mathcal{L}g(x,y)\le\rho[\delta_0x^{\rho_1}+y^{\rho_2}]^{\rho-1}
[-J_1(x,y)+J_2(x,y)]\le0$ for all $x,y>0$, which concludes the assertion.
\qed

\noindent{\it Proof of Theorem \ref{t1.111}.}
The assertion (iii) follows from
Lemma \ref{t8.10} and Proposition \ref{t3.1}.
The conclusion of (i)
is given by Lemma \ref{t5.4}(iii) and Proposition \ref{t3.1b}.
The conclusion of (ii) can be obtained
by Lemma \ref{t2.3} and Proposition \ref{t2.4b}.
\qed


\begin{thebibliography}{99}

		\bibitem{Bao}
		Bao, J.,  Mao, X., Yin, G. and Yuan C. (2011): Competitive Lotka-Volterra population dynamics with jumps. \textit{Nonlinear
			Anal.} {\bf 74}.  6601--6616.

		\bibitem{BaoShao}
Bao, J. and Shao, J. (2016): Permanence and extinction of regime-switching predator-prey models. \textit{SIAM J. Math. Anal.} \textbf{48}, 725--739.

\bibitem{BarnesMytnikSun}
Barnes, C., Mytnik, L., and Sun, Z. (2024):
On the coming down from infinity of coalescing Brownian motions.
\textit{Ann. Probab.}, \textbf{52}, 67--92.

\bibitem{Cattiaux}
Cattiaux, P. and M\'el\'eard, S. (2010): Competitive or weak cooperative stochastic Lotka-Volterra systems conditioned on
non-extinction. \textit{J. Math. Biol.} {\bf60}, 797--829.


\bibitem{Chen1986a}
Chen, M. (1986a): Couplings of jump processes. \textit{Acta Math. Sinica, New Series.} \textbf{2}, 121--136.


\bibitem{Chen1986b}
Chen, M. (1986b): \textit{Jump Processes and Interacting Particle Systems} (in Chinese). Beijing Normal Univ. Press.




\bibitem{Chen04}
Chen, M. (2004):
\textit{From Markov Chains to Non-Equilibrium Particle Systems}.
Second ed., World Scientific.


\bibitem{ChenYangZhou}
Chen, S., Yang, X. and Zhou, X. (2025):
Extinction, explosion and contraction for time-inhomogeneous SDEs with jumps.
ArXiv:2511.16973.



		\bibitem{DuDangYin} Du, N,  Dang, N. and Yin, G. (2016):
		Conditions for permanence and ergodicity of certain
		stochastic predator-prey models. \textit{J. Appl. Probab.}  \textbf{53}, 187--202.



\bibitem{Evans}
Evans, S. N., Hening, A. and Schreiber, S. (2015):
Protected polymorphisms and evolutionary stability of patch-selection
strategies in stochastic environents. \textit{J. Math. Biol.} {\bf71}, 325--359.



\bibitem{FoucartLiZhou21}
Foucart, C. Li, P.-S., and Zhou, X. (2021). Time-changed spectrally positive L\'evy processes started from infinity. \textit{Bernoulli}, \textbf{27(2)}, 1291--1318.

\bibitem{Grey74}
		Grey, D. R. (1974): Asymptotic behaviour of continuous time, continuous state-space branching processes. \textit{J. Appl. Probab.} \textbf{11}, 669--677.

\bibitem{Ky}
Kyprianou, A. (2006): \textit{Introductory Lectures on Fluctuations of L\'evy
Processes with Applications.} Universitext. Springer-Verlag, Berlin.

\bibitem{LiZh21}
Li, B. and Zhou, X. (2021):
On the explosion of a class of continuous-state
nonlinear branching processes. \textit{Electron. J. Probab.} \textbf{26}, no. 148, 1--25.



\bibitem{LiP19}
		Li, P.-S.  (2019): A continuous-state polynomial branching process. \textit{Stochastic Process. Appl.} \textbf{129}, 2941--2967.

\bibitem{LYZh}
Li, P.-S., Yang, X. and Zhou, X. (2019):
A general continuous-state nonlinear branching process.
\textit{Ann. Appl. Probab.} \textbf{29}, 2523-2555.




\bibitem{Li2020}
Li, Z. (2020): \textit{Continuous-state branching processes with immigration}. A Chapter in: From Probability to Finance, pp.1-69, edited by Y. Jiao. Mathematical Lectures from Peking University. Springer, Singapore.


\bibitem{Ma2013}
Ma, R. (2013): Stochastic equations for two-type continuous-state branching processes with immigration. \textit{Acta Math. Sin. (Engl. Ser.)}. {\bf29}, 287--294.

\bibitem{Pardoux16}
Pardoux, E. (2016):
\textit{Probabilistic Models of Population Evolution:
Scaling Limits, Genealogies and Interactions}. Springer, Heidelberg.

\bibitem{RXYZ19}
Ren, Y., Xiong, J., Yang, X. and Zhou, X. (2022):
On the extinction-extinguishing dichotomy for a stochastic
Lotka-Volterra type population dynamical system.
\textit{Stochastic Process. Appl.} \textbf{150}, 50-90.

\bibitem{Watanabe}
Watanabe, S. (1969):
On two dimensional Markov processes with branching property.
\textit{Trans. Amer. Math. Soc}. \textbf{136}, 447--466.

\bibitem{XYZh24}
Xiong, J., Yang, X. and Zhou, X. (2026):
Extinction behaviour for competing continuous-state  population dynamics.
\textit{Stochastic Process. Appl.} \textbf{196}, 104928.

\bibitem{XYZh25}
Xiong, J., Yang, X. and Zhou, X. (2026+): Extinction behaviour for  mutually
enhancing continuous-state population dynamics.
Submitted. ArXiv: 2603.16190v2.

\end{thebibliography}
\end{document}